\documentclass{svmult}

\usepackage{makeidx}         
\usepackage{graphicx}        
\usepackage{multicol}        

\makeindex             

\usepackage{amssymb}
\usepackage{amsmath}
\usepackage{xypic}

\newcommand{\C}{\ensuremath{\mathbb{C}}}
\newcommand{\R}{\ensuremath{\mathbb{R}}}
\newcommand{\Q}{\ensuremath{\mathbb{Q}}}
\newcommand{\Z}{\ensuremath{\mathbb{Z}}}

\newcommand{\cale}{{\cal E}}
\newcommand{\calf}{{\cal F}}

\newcommand{\calh}{{\cal H}}

\newcommand{\calk}{{\cal K}}
\newcommand{\call}{{\cal L}}
\newcommand{\calm}{{\cal M}}
\newcommand{\caln}{{\cal N}}

\newcommand{\calu}{{\cal U}}

\newcommand{\norm}[1]{\ensuremath{\|#1\|}}
\newcommand{\inner}[1]{\ensuremath{\langle #1 \rangle}}
\newcommand{\csa}{$C^*$-algebra}
\newcommand{\subcsa}{sub-$C^*$-algebra}
\newcommand{\Gcsa}{$G$-$C^*$-algebra}
\newcommand{\EG}{\ensuremath{\underline{E}G}}

\newcounter{mnote}
\newenvironment{metanote}[1][Note]{\begin{quote}\stepcounter{mnote}
\textbf{\textsc{Note}}~\sl\marginpar{\textsc{#1}
\arabic{mnote}}}{\end{quote}}

\begin{document}

\title*{$K$-theory for group \csa s}

\author{Paul F.~Baum\inst{1} \and Rub\'en J.~S\'anchez-Garc\'ia\inst{2}}

\institute{Mathematics Department\\ 206 McAllister Building\\ The
Pennsylvania State University\\ University Park, PA 16802\\ USA\\
\texttt{baum@math.psu.edu}
\and
Mathematisches Institut\\Heinrich-Heine-Universit\"at D\"usseldorf\\
Universit\"atsstr. 1\\
40225 D\"usseldorf\\Germany\\
\texttt{sanchez@math.uni-duesseldorf.de}}

%
%
\maketitle


\section*{Introduction}
These notes are based on a lecture course given by the first author in the Sedano Winter School on $K$-theory held in Sedano, Spain, on January 22-27th of 2007. They aim at introducing $K$-theory of \csa{}s, equivariant $K$-homology and $KK$-theory in the context of the Baum-Connes conjecture.

We start by giving the main definitions, examples and properties of \csa{}s in Section \ref{section:Csas}. A central construction is the reduced \csa{} of a locally compact, Hausdorff, second countable group $G$. In Section \ref{section:KTheoryofcsas} we define $K$-theory for \csa{}s, state the Bott periodicity theorem and establish the connection with Atiyah-Hirzebruch topological $K$-theory.

Our main motivation will be to study the $K$-theory of the reduced \csa{} of a group $G$ as above. The Baum-Connes conjecture asserts that these $K$-theory groups are isomorphic to the equivariant $K$-homology groups of a certain $G$-space, by means of the index map. The $G$-space is the universal example for proper actions of $G$, written $\EG$. Hence we procceed by discussing proper actions in Section \ref{section:ProperSpaces} and the universal space \EG{} in Section \ref{section:ClassifyingSpace}. 

Equivariant $K$-homology is explained in Section \ref{section:EquivariantKHomology}. This is an equivariant version of the dual of Atiyah-Hirzebruch $K$-theory. Explicitly, we define the groups $K^G_j(X)$ for $j = 0,1$ and $X$ a proper $G$-space with compact, second countable quotient $G\backslash X$. These are quotients of certain equivariant $K$-cycles by homotopy, although the precise definition of homotopy is postponed. We then address the problem of extending the definition to \EG{}, whose quotient by the $G$-action may not be compact.

In Section \ref{section:DiscreteCase} we concentrate on the case when $G$ is a discrete group, and in Section \ref{section:CompactCase} on the case $G$ compact. In Section \ref{section:EquivariantKhomology2} we introduce $KK$-theory for the first time. This theory, due to Kasparov, is a generalization of both $K$-theory of \csa{}s and $K$-homology. Here we define $KK^j_G(A,\mathbb{C})$ for a separable \csa{} $A$ and $j=0,1$, although we again postpone the exact definition of homotopy. The already defined $K^G_j(X)$ coincides with this group when $A=C_0(X)$.

At this point we introduce a generalization of the conjecture called the Baum-Connes conjecture with coefficients, which consists in adding coefficients in a \Gcsa{} (Section \ref{section:BCCwC}). To fully describe the generalized conjecture we need to introduce Hilbert modules and the reduced crossed-product (Section \ref{section:HilbertModules}), and to define $KK$-theory for pairs of \csa{}s. This is done in the non-equivariant situation in Section \ref{section:HomotopyAndKKtheory} and in the equivariant setting in Section \ref{section:EquivariantKhomology3}. In addition we give at this point the missing definition of homotopy. Finally, using equivariant $KK$-theory, we can insert coefficients in equivariant $K$-homology, and then extend it again to \EG.

The only ingredient of the conjecture not yet accounted for is the index map. It is defined in Section \ref{section:IndexMap} via the Kasparov product and descent maps in $KK$-theory. We finish with a brief exposition of the history of $K$-theory and a discussion of Karoubi's conjecture, which symbolizes the unity of $K$-theory, in Section \ref{section:History}.

\section{\csa s}\label{section:Csas}
\label{sec:CstarAlgebras} We start with some definitions and basic
properties of \csa{}s. Good references for \csa{} theory are \cite{Arveson76}, \cite{Dixmier77},
\cite{Murphy90} or \cite{Pedersen79}.

\subsection{Definitions}
\begin{definition}
A \emph{Banach algebra} is an (associative, not necessarily unital)
algebra $A$ over \C{} with a given norm $\|~\|$
\[
    \|~\| : A \longrightarrow [0,\infty)
\]
such that $A$ is a complete normed algebra, that is, for all $a,b \in A$, $\lambda \in \C$,
\begin{enumerate}
  \item $\norm{\lambda a} = |\lambda| \norm{a}$,
  \item $\norm{a+b} \le \norm{a} + \norm{b}$,
  \item $\norm{a} = 0 \Leftrightarrow a = 0$,
  \item $\norm{ab} \le \norm{a}\norm{b}$,
  \item every Cauchy sequence is convergent in $A$ (with respect to
  the metric $d(a,b)=\norm{a-b}$).
\end{enumerate}
\end{definition}


A \csa{} is a Banach algebra with an
involution satisfying the \emph{\csa{} identity}.
\begin{definition}
A \emph{\csa{}} $A = (A,\norm{~}, *)$ is a Banach algebra $(A,
\norm{~})$ with a map $* : A \rightarrow A, a \mapsto a^*$ such that for all $a,b \in A$, $\lambda \in \C$
\begin{enumerate}
  \item $(a+b)^* = a^*+b^*$,
  \item $(\lambda a)^* = \overline{\lambda} a^*$,
  \item $(ab)^* = b^*a^*$,
  \item $(a^*)^* = a$,
  \item $\norm{aa^*} = \norm{a}^2$ (\emph{\csa{} identity}).
\end{enumerate}
\end{definition}
Note that in particular $\norm{a} = \norm{a^*}$ for all $a \in A$:
 for $a=0$ this is clear; if $a \neq 0$ then $\|a\| \neq 0$
and $\|a\|^2 = \| a a^* \| \le \|a\| \|a^*\|$ implies $\|a\| \le
\|a^*\|$, and similarly $\|a^*\| \le \|a\|$.

A \csa{} is \emph{unital} if it has a multiplicative unit $1 \in A$.
A \emph{sub-$C^*$-algebra} is a non-empty subset of $A$ which is a
$C^*$-algebra with the operations and norm given on $A$.
\begin{definition}
A \emph{$*$-homomorphism} is an algebra homomorphism $\varphi: A
\rightarrow B$ such that $\varphi(a^*) = (\varphi(a))^*$, for all $a
\in A$.
\end{definition}
\begin{proposition}
If $\varphi: A \rightarrow B$ is a $*$-homomorphism then
$\norm{\varphi(a)} \le \norm{a}$ for all $a \in A$. In particular,
$\varphi$ is a (uniformly) continuous map.
\end{proposition}
For a proof see, for instance, \cite[Thm. 1.5.7]{Pedersen79}.


\subsection{Examples}
We give three examples of $C^*$-algebras.
\begin{example}\label{example:C_0(X)}
Let $X$ be a Hausdorff, locally compact topological space. Let $X^+
= X \cup \{p_\infty\}$ be its one-point compactification. (Recall
that $X^+$ is Hausdorff if and only if $X$ is Hausdorff and locally
compact.)

Define the $C^*$-algebra
\[
    C_0(X) = \left\{ \alpha : X^+ \rightarrow \C \, | \, \alpha
    \textrm{ continuous, } \alpha(p_\infty) = 0 \right\}\,,
\]
with operations: for all $\alpha, \beta \in C_0(X), p \in X^+,
\lambda \in \C$
\begin{eqnarray*}
  (\alpha + \beta)(p) &=& \alpha(p) + \beta(p), \\
  (\lambda\alpha)(p) &=& \lambda\alpha(p),\\
  (\alpha\beta)(p) &=& \alpha(p)\beta(p), \\
  \alpha^*(p) &=& \overline{\alpha(p)},\\
  \|\alpha\| &=& \sup_{p \in X} |\alpha(p)|\,.
\end{eqnarray*}
Note that if $X$ is compact Hausdorff, then
\[
    C_0(X) = C(X) = \{ \alpha : X \rightarrow \C \,|\, \alpha
    \textrm{ continuous}\}\,.
\]
\end{example}


\begin{example}
Let $H$ be a Hilbert space. A Hilbert space is \emph{separable} if
it admits a countable (or finite) orthonormal basis. (We shall deal with
separable Hilbert spaces unless explicit mention is made to the contrary.)

Let $\call(H)$ be the set of bounded linear operators on $H$, that
is, linear maps $T : H \rightarrow H$ such that
\[
    \norm{T} = \sup_{\substack{\norm{u}=1}} \norm{Tu} < \infty\,,
\]
where $\norm{u} = \langle u, u \rangle^{1/2}$. It is a complex
algebra with
\begin{eqnarray*}
   (T+S)u &=& T u + S u,\\
   (\lambda T) u &=& \lambda(Tu),\\
   (TS)u &=& T(Su),
\end{eqnarray*}
for all $T,S \in \call(H)$, $u \in H$, $\lambda \in \C$. The norm is
the operator norm $\|T\|$ defined above, and $T^*$ is the
adjoint operator of $T$, that is, the unique bounded operator such
that
\[
    \inner{Tu, v} = \inner{u, T^*v}
\]
for all $u,v \in H$.
\end{example}


\begin{example}
Let $\call(H)$ be as above. A bounded operator is \emph{compact} if
it is a norm limit of operators with finite-dimensional image, that is,
\[
    \calk(H) = \{ T \in \call(H) \,|\, T \textup{ compact
    operator}\} = \overline{\{ T \in \call(H) \,|\, \dim_\C T(H) < \infty
    \}}\,,
\]
where the overline denotes closure with respect to the operator norm.
$\calk(H)$ is a \subcsa{} of $\call(H)$. Moreover, it is an
ideal of $\call(H)$ and, in fact, the only norm-closed ideal except
0 and $\call(H)$.
\end{example}


\subsection{The reduced \csa{} of a group}\label{section:ReducedCsaGroup} Let $G$ be a
topological group which is locally compact, Hausdorff and second
countable (i.e. as a topological space it has a countable basis). There is a
\csa{} associated to $G$, called the \emph{reduced \csa{} of $G$}, defined as follows.

\begin{remark} We need $G$ to be locally compact and Hausdorff to
guarantee the existence of a Haar measure. The countability
assumption makes the Hilbert space $L^2(G)$ separable and also
avoids some technical difficulties when later defining Kasparov's
$KK$-theory.
\end{remark}

Fix a left-invariant Haar measure $dg$ on $G$. By left-invariant we mean that if $f\colon G
\rightarrow \C$ is continuous with compact support then
\[
    \int_G f(\gamma g) dg = \int_G f(g) dg \qquad \textrm{for all } \gamma \in
    G\,.
\]
Define the Hilbert space $L^2G$ as
\begin{eqnarray*}
  L^2G &=& \Big\{ u:G\rightarrow \C \ \big|\, \int_G |u(g)|^2 dg < \infty
    \Big\},
\end{eqnarray*}
with scalar product
\[
    \inner{u,v} = \int_G \overline{u(g)}\,v(g)dg
\]
for all $ u,v \in L^2G$.


Let $\call(L^2G)$ be the \csa{} of all bounded linear operators $T: L^2G
\rightarrow L^2G$. On the other hand, define
\[
    C_cG = \left\{ f:G\rightarrow \C \,|\, f
    \textup{ continuous with compact support} \right\}.
\]
It is an algebra with
\begin{eqnarray*}
  (f + h)(g) &=& f(g) + h(g),\\
  (\lambda f)(g) &=& \lambda f(g),
\end{eqnarray*}
for all $f,h \in C_cG$, $\lambda \in \C$, $g \in G$, and
multiplication given by \emph{convolution}
\[
    (f \ast h) (g_0) = \int_G f(g)h(g^{-1}g_0)\,dg \quad \textup{for all } g_0 \in G.
\]
\begin{remark}
  When $G$ is discrete, $\int_G f(g) dg = \sum_G f(g)$ is a Haar measure,  $C_cG$ is the complex group algebra $\mathbb{C}[G]$ and $f \ast h$ is the usual product in $\mathbb{C}[G]$.
\end{remark}

There is an injection of algebras
\[\begin{array}{ccccc}
    0 &\longrightarrow& C_cG &\longrightarrow& \call(L^2G)\\
     & & f &\mapsto& T_f
\end{array}\]
where
\begin{align*}
    &T_f(u) = f \ast u &u \in L^2G\,,\\
    &(f \ast u) (g_0) = \int_G f(g)u(g^{-1}g_0)dg &g_0 \in G\,.
\end{align*}
Note that $C_cG$ is not necessarily a sub-\csa{} of
$\call(L^2G)$ since it may not be complete. We define $C^*_r(G)$,
the \emph{reduced \csa{} of $G$}, as the norm closure of $C_cG$ in
$\call(L^2G)$:
\[
    C^*_r(G) = \overline{C_cG} \subset \call(L^2G).
\]
\begin{remark}
There are other possible completions of $C_cG$. This particular one, i.e. $C^*_r(G)$, uses only the
left regular representation of $G$ (cf. \cite[Chapter 7]{Pedersen79}).
\end{remark}

\subsection{Two classical theorems}
We recall two classical theorems about \csa{}s. The first one says that any
\csa{} is (non-canonically) isomorphic to a \csa{} of operators, in the
sense of the following definition.
\begin{definition}
    A subalgebra $A$ of $\call(H)$ is a \emph{\csa{} of
    operators} if
    \begin{enumerate}
      \item $A$ is closed with respect to the operator norm;
      \item if $T \in A$ then the adjoint operator $T^* \in A$.
    \end{enumerate}
\end{definition}
That is, $A$ is a \subcsa{} of $\call(H)$, for some Hilbert
space  $H$.

\begin{theorem}[I.~Gelfand and V.~Naimark]
Any \csa{} is isomorphic, as a \csa{}, to a \csa{} of operators.
\end{theorem}

The second result states that any commutative \csa{} is
(canonically) isomorphic to $C_0(X)$,
for some topological space $X$.
\begin{theorem}[I.~Gelfand]\label{thm:Gelfand}
    Let $A$ be a commutative \csa{}. Then $A$ is (canonically)
    isomorphic to $C_0(X)$ for $X$ the space of maximal ideals of
    $A$.
\end{theorem}

\begin{remark} The topology on $X$ is the \emph{Jacobson topology}
or \emph{hull-kernel topology} \cite[p.~159]{Murphy90}.
\end{remark}
Thus a non-commutative \csa{} can be viewed as a `non-commutative,
locally compact, Hausdorff topological space'.

\subsection{The categorical viewpoint}
Example \ref{example:C_0(X)} gives a functor between the category of locally compact,
Hausdorff, topological spaces and the category of \csa s, given by $X \mapsto C_0(X)$.
Theorem \ref{thm:Gelfand} tells us that its restriction to commutative \csa
s is an equivalence of categories,
\begin{align*}
    \left( \begin{array}{c} \textrm{commutative}\\
                            \textrm{\csa s}\\
                            \end{array}
    \right)
    & \simeq
        \left( \begin{array}{c} \textrm{locally compact, Hausdorff,}\\
                            \textrm{topological spaces}\\
                            \end{array}
    \right)^{op} \\
    C_0(X) & \longleftarrow  X
\end{align*}

On one side we have \csa s and $*$-homorphisms, and on the other
locally compact, Hausdorff topological spaces with morphisms from $Y$ to $X$ being continuous maps $f \colon X^+ \rightarrow Y^+$ such that $f(p_\infty)=q_\infty$. (The symbol $op$ means the opposite or dual category, in other words, the functor is contravariant.)

\begin{remark}
This is not the same as continuous proper maps $f \colon X \rightarrow Y$ since we do not require that the map $f \colon X^+ \rightarrow Y^+$ maps $X$ to $Y$.
\end{remark}


\section{$K$-theory of \csa s}\label{section:KTheoryofcsas}
In this section we define the $K$-theory groups of an arbitrary \csa. We first give the definition for a \csa{} with unit and then extend it to the non-unital case. We also discuss Bott periodicity and the connection with topological $K$-theory of spaces. More details on $K$-theory of \csa{}s is given in Section 3 of Corti\~nas' notes \cite{Willie}, including a proof of Bott periodicity. Other references are \cite{Murphy90}, \cite{RLL00} and \cite{Wegge-Olsen93}.

Our main motivation is to study the $K$-theory of $C_r^*(G)$, the
reduced \csa{} of $G$. From Bott periodicity, it suffices to compute $K_j\left( C_r^*(G) \right)$ for $j=0,1$.
In 1980, Paul Baum and Alain Connes conjectured that
these $K$-theory groups are isomorphic to the \emph{equivariant
$K$-homology} (Section \ref{section:EquivariantKHomology}) of a
certain $G$-space. This $G$-space is the \emph{universal example for
proper actions of $G$} (Sections \ref{section:ProperSpaces} and \ref{section:ClassifyingSpace}), written
\EG{}. Moreover, the conjecture states that the isomorphism is given
by a particular map called the \emph{index map} (Section \ref{section:IndexMap}).

\begin{conjecture}[P.~Baum and A.~Connes, 1980]\label{conj:BCC}
Let $G$ be a locally compact, Hausdorff, second countable,
topological group. Then the index map
\[
    \mu : K^G_j\left(\EG\right) \longrightarrow K_j\left( C_r^*(G)
    \right) \quad j=0,1
\]
is an isomorphism.
\end{conjecture}


\subsection{Definition for unital \csa{}s}
Let $A$ be a \csa{} with unit $1_A$. Consider $GL(n,A)$, the group
of invertible $n$ by $n$ matrices with coefficients in $A$. It is a
topological group, with topology inherited from $A$. We have a
standard inclusion
\begin{eqnarray*}
  GL(n,A) &\hookrightarrow& GL(n+1,A) \\
    \begin{pmatrix}
        a_{11} & \ldots & a_{1n}\\
        \vdots & \cdots & \vdots\\
        a_{n1} & \ldots & a_{nn}
    \end{pmatrix}
        &\mapsto&
    \begin{pmatrix}
        a_{11} & \ldots & a_{1n} & 0\\
        \vdots & \cdots & \vdots & \vdots\\
        a_{n1} & \ldots & a_{nn} & 0\\
        0 & \ldots & 0 & 1_A
    \end{pmatrix}\,.
\end{eqnarray*}
Define $GL(A)$ as the direct limit with respect to these inclusions
\[
    GL(A) = \bigcup_{n=1}^\infty GL(n,A)\,.
\]
It is a topological group with the \emph{direct limit topology}: a subset $\theta$ is open if and only
if $\theta \cap GL(n,A)$ is open for every $n \ge 1$. In
particular, $GL(A)$ is a topological space, and hence we can
consider its homotopy groups.
\begin{definition}[$K$-theory of a unital \csa]
$$K_j(A) = \pi_{j-1}\left( GL(A) \right) \quad j = 1, 2, 3, \ldots$$
\end{definition}
Finally, we define $K_0(A)$ as the \emph{algebraic} $K$-theory group
of the ring $A$, that is, the Grothendieck
group of finitely generated (left) projective $A$-modules (cf. \cite[Remark 2.1.9]{Willie}),
$$K_0(A) = K_0^{\textrm{alg}}(A)\,.$$
\begin{remark} Note that $K_0(A)$ only depends on the ring structure
of $A$ and so we can `forget' the norm and the involution.
The definition of $K_1(A)$ does require the norm but not
the involution, so in fact we are defining $K$-theory of Banach algebras with unit. Everything we say in \ref{subsect:BottPeriodicity} below, including Bott periodicity, is true for Banach algebras.
\end{remark}

\subsection{Bott periodicity}\label{subsect:BottPeriodicity}
The fundamental result is \emph{Bott
periodicity}. It says that the homotopy groups of $GL(A)$ are
periodic modulo 2 or, more precisely, that the double loop space of
$GL(A)$ is homotopy equivalent to itself,
\[
    \Omega^2 GL(A) \simeq GL(A)\,.
\]
As a consequence, the $K$-theory of the \csa{} $A$ is periodic modulo 2
\[
    K_j(A) = K_{j+2}(A) \quad j \ge 0.
\]
Hence from now on we will only consider $K_0(A)$ and $K_1(A)$.

\subsection{Definition for non-unital \csa s} If $A$ is a \csa{} without a unit,
we formally adjoin one.
Define $\widetilde{A} = A \oplus \C$ as a complex algebra with
multiplication, involution and norm given by
\begin{eqnarray*}
    (a,\lambda) \cdot (b, \mu) &=& (ab + \mu a + \lambda b, \lambda
    \mu), \\
    (a,\lambda)^* &=& (a^*,\overline{\lambda}),\\
    \norm{(a,\lambda)} &=& \sup_{\norm{b} = 1} \norm{a b + \lambda b}\,.
\end{eqnarray*}
This makes $\widetilde{A}$ a unital \csa{} with unit $(0,1)$.
We have an exact sequence
\[
    0 \longrightarrow A \longrightarrow \widetilde{A}
    \longrightarrow \C \longrightarrow 0.
\]

\begin{definition} Let $A$ be a non-unital \csa. Define $K_0(A)$ and $K_1(A)$
as
\begin{eqnarray*}
  K_0(A) &=& \ker \left( K_0(\widetilde{A}) \rightarrow K_0(\C)
    \right) \\
  K_1(A) &=& K_1(\widetilde{A}).
\end{eqnarray*}
\end{definition}
This definition agrees with the previous one when $A$ has a unit. It also satisfies Bott periodicity (see Corti\~nas' notes \cite[3.2]{Willie}).

\begin{remark} Note that the \csa{} $C^*_r(G)$ is unital if and only if $G$ is
discrete, with unit the Dirac function on $1_G$.
\end{remark}

\begin{remark}
There is algebraic $K$-theory of rings (see \cite{Willie}). Althought a \csa{} is in particular a ring, the two $K$-theories are different; algebraic $K$-theory does not satisfy Bott periodicity and $K_1$ is in general a quotient of $K_1^{alg}$. We shall compare both definitions in Section \ref{section:Comparison} (see also \cite[Section 7]{Willie}).
\end{remark}



\subsection{Functoriality}
Let $A, B$ be \csa s (with or without units), and $\varphi: A
\rightarrow B$ a $*$-homomorphism. Then $\varphi$ induces a
 homomorphism of abelian groups
\[
    \varphi_* : K_j(A) \longrightarrow K_j(B) \quad j=0,1.
\]
This makes $A \mapsto K_j(A)$, $j=0,1$,  covariant functors from \csa s to
abelian groups \cite[Sections 4.1 and 8.2]{RLL00}.

\begin{remark}
When $A$ and $B$ are unital and $\varphi(1_A)=1_B$, the map $\varphi_*$ is the one induced by $GL(A) \rightarrow GL(B)$, $(a_{ij}) \mapsto \left(\varphi(a_{ij})\right)$ on homotopy groups.
\end{remark}


\subsection{More on Bott periodicity}
In the original article \cite{Bott59}, Bott computed the stable homotopy of the classical groups and, in particular, the homotopy groups $\pi_j(GL(n,\C)$ when $n \gg j$.

\begin{figure}[h]\begin{center}\includegraphics[width=4.5cm, height=6cm]{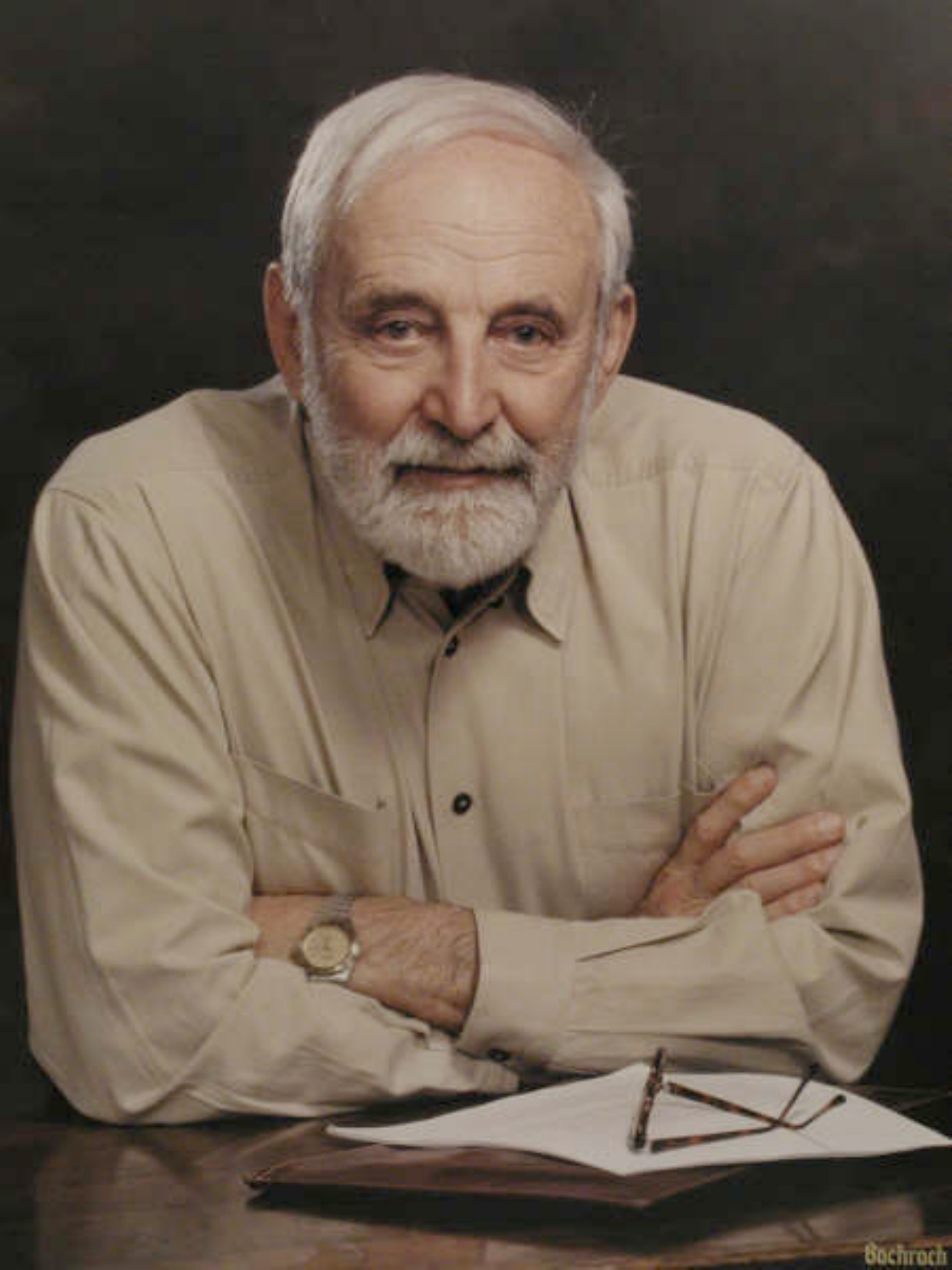}\caption{Raoul Bott}\end{center}\end{figure}

\begin{theorem}[R.~Bott \cite{Bott59}] The homotopy groups of $GL(n,\C)$ are
\[
    \pi_j\left( GL(n,\C) \right) = \begin{cases} 0 &
    \textrm{$j$ even}\\ \Z & \textrm{$j$ odd} \end{cases}
\]
for all $j=0,1,2, \ldots, 2n-1$.
\end{theorem}


As a corollary of the previous theorem, we obtain the $K$-theory of \C,
considered as a \csa.
\begin{theorem}[R.~Bott]\label{Thm:Bott2}
    $$K_j(\C) = \begin{cases} \Z & \textrm{$j$ even,}\\ 0 & \textrm{$j$ odd.}
    \end{cases}$$
\end{theorem}
\textit{Sketch of proof.}
Since \C{} is a field, $K_0(\C) = K_0^{alg}(\C) = \Z$. By the polar decomposition, $GL(n,\C)$ is homotopy equivalent to $U(n)$. The homotopy long exact sequence of the fibration $U(n) \rightarrow U(n+1) \rightarrow S^{2n+1}$ gives $\pi_j(U(n)) = \pi_j(U(n+1))$ for all $j \le 2n+1$. Hence $K_j(\C)= \pi_{j-1}(GL(\C)) = \pi_{j-1}(GL(2j-1,\C))$ and apply the previous theorem.


\begin{remark}
Compare this result with $K_1^{alg}(\C) = \C^\ast$ (since \C{} is a field, see \cite[Ex. 3.1.6]{Willie}). Higher algebraic $K$-theory groups for $\mathbb{C}$ are only partially understood.
\end{remark}


\subsection{Topological $K$-theory}\label{section:TopologicalKTheory}
There is a close connection between $K$-theory
of \csa s and topological $K$-theory of spaces.

Let $X$ be a locally compact, Hausdorff, topological space. Atiyah and
Hirzebruch \cite{AH59} defined abelian groups $K^0(X)$ and $K^1(X)$ called
\emph{topological $K$-theory with compact supports}. For instance, if $X$ is compact,
$K^0(X)$ is the Grothendieck group of complex vector bundles on $X$.

\begin{theorem} Let $X$ be a locally compact, Hausdorff, topological
space. Then
\[
    K^j(X) = K_j\left(C_0(X)\right), \ j=0,1.
\]
\end{theorem}
\begin{remark} This is known as Swan's theorem when $j=0$ and $X$ compact.
\end{remark}
In turn, topological $K$-theory can be computed up to torsion via a Chern character.
Let $X$ be as above. There is a \emph{Chern character} from
topological $K$-theory to rational cohomology with compact supports
\[
    ch : K^j(X) \longrightarrow \bigoplus_{l\ge 0} H_c^{j+2l} (X;\Q)\,, \quad
    j=0,1.
\]

Here the target cohomology theory $H^*_c(-;\Q)$ can be \v{C}ech cohomology
with compact supports, Alexander-Spanier cohomology with compact supports
or representable Eilenberg-MacLane cohomology with compact supports.

This map becomes an isomorphism when tensored with the rationals.
\begin{theorem} Let $X$ be a locally compact, Hausdorff, topological space.
The Chern character is a rational isomorphism, that is,
\[
    K^j(X)\otimes_\Z \Q \longrightarrow \bigoplus_{l\ge 0} H_c^{j+2l} (X;\Q)\,, \quad
    j=0,1
\]
is an isomorphism.
\end{theorem}
\begin{remark} This theorem is still true for singular cohomology when $X$ is a locally finite CW-complex.
\end{remark}



\section{Proper $G$-spaces}\label{section:ProperSpaces}
In the following three sections, we will describe
the left-hand side of the Baum-Connes conjecture (Conjecture
\ref{conj:BCC}). The space \EG{} appearing on the topological side
of the conjecture is the \emph{universal example for proper actions
for $G$}. Hence we will start by studying proper $G$-spaces.

Recall the definition of $G$-space, $G$-map and $G$-homotopy.
\begin{definition}
A \emph{$G$-space} is a topological space $X$ with a given
continuous action of $G$
\[
    G \times X \longrightarrow X.
\]
A \emph{$G$-map} is a continuous map $f : X \rightarrow Y$ between
$G$-spaces such that
\[
    f(gp) = gf(p) \ \textup{for all } (g,p) \in G \times X.
\]
Two $G$-maps $f_0, f_1 : X \rightarrow Y$ are
\emph{$G$-homotopic} if they are homotopic through $G$-maps, that
is, there exists a homotopy $\{f_t\}_{0 \le t \le 1}$ with each
$f_t$ a $G$-map.
\end{definition}

We will require proper $G$-spaces to be Hausdorff and paracompact.
Recall that a space $X$ is \emph{paracompact} if every open cover of
$X$ has a locally finite open refinement or, alternatively, a locally
finite partition of unity subordinate to any given open cover.

\begin{remark} Any metrizable space (i.e.~there is a metric with the
same underlying topology) or any CW-complex (in its usual CW-topology) is
Hausdorff and paracompact. 
\end{remark}


\begin{definition}
A $G$-space $X$ is \emph{proper} if
\begin{itemize}
  \item $X$ is Hausdorff and paracompact;
  \item the quotient space $G\backslash X$ (with the quotient topology) is
  Hausdorff and paracompact;
  \item for each $p \in X$ there exists a triple $(U,H,\rho)$ such
  that
  \begin{enumerate}
    \item $U$ is an open neighborhood of $p$ in $X$ with $gu \in U$
    for all $(g,u)\in G\times U$;
    \item $H$ is a compact subgroup of $G$;
    \item $\rho : U \rightarrow G/H$ is a $G$-map.
  \end{enumerate}
\end{itemize}
\end{definition}
Note that, in particular, the stabilizer $\textup{stab}(p)$ is a closed subgroup of a conjugate of $H$ and hence compact.

\begin{remark}
The converse is not true in general; the action of $\Z$ on $S^1$ by an irrational rotation is free but it is not a proper $\Z$-space.
\end{remark}

\begin{remark}
If $X$ is a \emph{$G$-CW-complex} then it is a proper $G$-space (even in the weaker definition below) if and only if all the cell stabilizers are compact, see Thm.~1.23 in \cite{Luck89}.
\end{remark}

Our definition is stronger than the usual definition of proper $G$-space, which requires the map $G \times X \rightarrow X \times X$, $(g,x) \mapsto (gx,x)$ to be proper, in the sense that the pre-image of a compact set is compact.
Nevertheless, both definitions agree for locally compact, Hausdorff, second countable $G$-spaces.


\begin{proposition}[J.~Chabert, S.~Echterhoff, R.~Meyer \cite{CEM01}]
If $X$ is a locally compact, Hausdorff, second countable $G$-space,
then $X$ is proper if and only if the map
\begin{align*}
    G \times X &\longrightarrow X \times X\\
    (g,x) & \longmapsto (gx,x)
\end{align*}
is proper.
\end{proposition}

\begin{remark}
For a more general comparison among these and other definitions of proper actions see \cite{Biller04}.
\end{remark}

\section{Classifying space for proper actions}\label{section:ClassifyingSpace}
Now we are ready for the definition of the space \EG{} appearing in the statement of the Baum-Connes Conjecture. Most of the material in this section is based on Sections 1 and 2 of \cite{BCH94}.

\begin{definition}
A \emph{universal example for proper actions of $G$}, denoted \EG{}, is a proper $G$-space such that:
\begin{itemize}
\item if $X$ is any proper $G$-space, then there exists a $G$-map $f: X \rightarrow \EG$ and any two $G$-maps from $X$ to \EG{} are  $G$-homotopic.
\end{itemize}
\end{definition}
\EG{} exists for every topological group $G$ \cite[Appendix 1]{BCH94} and it is unique up to $G$-homotopy, as follows. Suppose that \EG{} and $(\EG)'$ are both universal examples for proper actions of $G$. Then there exist $G$-maps
\begin{eqnarray*}
    f &:& \EG \longrightarrow (\EG)' \\ f' &:& (\EG)' \longrightarrow \EG
\end{eqnarray*}
and  $f' \circ f$ and $f \circ f'$ must be $G$-homotopic to the identity maps of \EG{} and (\EG)' respectively.


The following are equivalent axioms for a space $Y$ to be \EG{} \cite[Appendix 2]{BCH94}.
\begin{enumerate}
\item  $Y$ is a proper $G$-space.
\item If $H$ is any compact subgroup of $G$, then there exists $p \in Y$ with $hp = p$ for all $h \in H$.
\item Consider $Y \times Y$ as a $G$-space via $g(y_0,y_1) = (gy_0, gy_1)$, and the maps
\begin{gather*}
\rho_0, \rho_1 \colon Y \times Y \longrightarrow Y\\
\rho_0(y_0,y_1)=y_0\,, \quad \rho_1(y_0,y_1)=y_1\,.
\end{gather*}
Then $\rho_0$ and $\rho_1$ are $G$-homotopic.
\end{enumerate}

\begin{remark}
It is possible to define a universal space for any family of (closed) subgroups of $G$ closed under conjugation and finite intersections \cite{LuckSurvey}. Then \EG{} is the universal space for the family of compact subgroups of $G$.
\end{remark}

\begin{remark}
The space \EG{} can always be assumed to be a $G$-CW-complex. Then there is a homotopy characterization: a proper $G$-CW-complex $X$ is an \EG{} if and only if for each compact subgroup $H$ of $G$ the fixed point subcomplex $X^H$ is contractible (see \cite{LuckSurvey}).
\end{remark}


\subsection*{Examples}
\begin{enumerate}
\item If $G$ is compact, \EG{} is just a one-point space.

\item If $G$ is a Lie group with finitely many connected components then $\EG = G/H$, where $H$ is a maximal compact subgroup (i.e.~maximal among compact subgroups). 

\item If $G$ is a $p$-adic group then $\EG = \beta G$ the affine Bruhat-Tits building for $G$. For example, $\beta SL(2,\Q_p)$ is the $(p+1)$-regular tree, that is, the unique tree with exactly $p+1$ edges at each vertex (see Figure \ref{fig:RegularTree}) (cf. \cite{SerreTrees03}).

\begin{figure}[h]
\begin{center}
\includegraphics[scale=0.8]{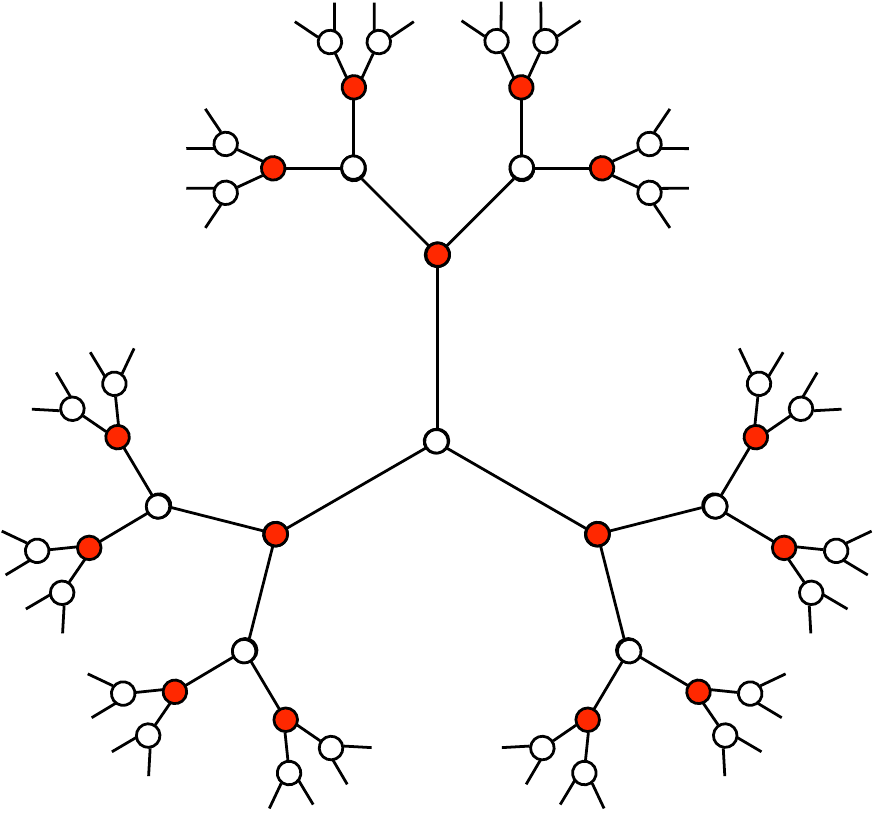}\label{fig:RegularTree}
\caption{The $(p+1)$-regular tree is $\beta SL(2,\Q_p)$}
\end{center}
\end{figure}

\item If $\Gamma$ is an arbitrary (countable) discrete group, there is an explicit construction,
$$
    \underline{E}\Gamma = \Big\{ f \colon \Gamma \rightarrow [0,1] \,\big|\, f \textup{ finite support }, \sum_{\gamma \in \Gamma} f(\gamma) = 1  \Big\}\,,
$$
that is, the space of all finite probability measures on $\Gamma$, topologized by the metric
$d(f,h) = \sqrt{\sum_{\gamma \in \Gamma} |f(\gamma) - h(\gamma)|^2}$.


\end{enumerate}



\section{Equivariant $K$-homology}\label{section:EquivariantKHomology}
$K$-homology is the dual theory to Atiyah-Hirzebruch $K$-theory (Section \ref{section:TopologicalKTheory}). Here we define an equivariant generalization due to Kasparov \cite{Kasparov80, Kasparov88}. If $X$ is a proper $G$-space with compact, second countable quotient then $K_i^G(X)$, $i=0,1$, are abelian groups defined as homotopy classes of $K$-cycles for $X$. These $K$-cycles can be viewed as $G$-equivariant abstract elliptic operators on $X$.

\begin{remark}
For a discrete group $G$, there is a topological definition of equivariant $K$-homology and the index map via equivariant spectra \cite{DL98}. This and other constructions of the index map are shown to be equivalent in \cite{HP04}.
\end{remark}

\subsection{Definitions}
Let $G$ be a locally compact, Hausdorff, second countable, topological group.

Let $H$ be a separable Hilbert space. Write $\mathcal{U}(H)$ for the set of unitary operators
$$
    \mathcal{U}(H) = \{ U \in \mathcal{L}(H) \, | \, U U^* = U^* U = I \}\,.
$$
\begin{definition}
    A \emph{unitary representation} of $G$ on $H$ is a group homomorphism $\pi \colon G \rightarrow \calu(H)$ such that for each $v \in H$ the map $\pi_v \colon G \rightarrow H$, $g \mapsto \pi(g)v$ is a continuous map from $G$ to $H$.
\end{definition}

\begin{definition}
A $G$-$C^*$-algebra is a \csa{} $A$ with a given continuous action of $G$
$$
    G \times A \longrightarrow A
$$
such that $G$ acts by $C^*$-algebra automorphisms.
\end{definition}
The continuity condition is that, for each $a \in A$, the map $G \rightarrow A$, $g \mapsto ga$ is a continuous map. We also have that, for each $g \in G$, the map $A \rightarrow A$, $a \mapsto ga$ is a \csa{} automorphism.


\begin{example} Let $X$ be a locally compact, Hausdorff $G$-space. The action of $G$ on $X$ gives an action of $G$ on $C_0(X)$,
$$
    (g\alpha)(x) = \alpha (g^{-1} x),
$$
where $g \in G$, $\alpha \in C_0(X)$ and $x \in X$. This action makes $C_0(X)$ into a $G$-\csa.
\end{example}
Recall that a \csa{} is \emph{separable} if it has a countable
dense subset.
\begin{definition}
Let $A$ be a separable \Gcsa. A \emph{representation} of $A$ is a triple $(H, \psi, \pi)$ with:
\begin{itemize}
    \item $H$ is a separable Hilbert space,
    \item $\psi \colon A \rightarrow \mathcal{L}(H)$ is a $*$-homomorphism,
     \item $\pi \colon G \rightarrow \mathcal{U}(H)$ is a unitary representation of $G$ on $H$,
     \item $\psi(ga) = \pi(g)\psi(a)\pi(g^{-1})$ for all $(g,a) \in G \times A$.
\end{itemize}
\end{definition}
\begin{remark}
We are using a slightly non-standard notation; in the literature this is usually called a \emph{covariant representation}.
\end{remark}

\begin{definition}
Let $X$ be a proper $G$-space with compact, second countable quotient space $G\backslash X$.
An \emph{equivariant odd $K$-cycle} for $X$ is a 4-tuple $(H,\psi, \pi, T)$ such that:
\begin{itemize}
    \item $(H,\psi,\pi)$ is a representation of the $G$-\csa{} $C_0(X)$,
    \item $T \in \mathcal{L}(H)$,
    \item $T = T^*$,
    \item $\pi(g)T - T\pi(g) = 0$ for all $g \in G$,
    \item $\psi(\alpha)T - T\psi(\alpha) \in \mathcal{K}(H)$ for all $\alpha \in C_0(X)$,
    \item $\psi(\alpha)(I - T^2) \in \mathcal{K}(H)$ for all $\alpha \in C_0(X)$.
\end{itemize}
\end{definition}
\begin{remark}
If $G$ is a locally compact, Hausdorff, second countable topological group and $X$ a proper $G$-space with locally compact quotient then $X$ is also locally compact and hence $C_0(X)$ is well-defined.
\end{remark}
Write $\cale_1^G(X)$ for the set of equivariant odd $K$-cycles for $X$. This concept was introduced by Kasparov as an abstraction an equivariant self-adjoint elliptic operator and goes back to Atiyah's theory of elliptic operators \cite{Atiyah70}.


\begin{example}
Let $G=\Z$, $X=\R$ with the action $\Z \times \R \rightarrow \R$, $(n,t) \mapsto n+t$. The quotient space is $S^1$, which is compact. Consider $H = L^2(\R)$ the Hilbert space of complex-valued square integrable functions with the usual Lebesgue measure. Let $\psi \colon C_0(\R) \rightarrow \call(L^2(\R))$ be defined as $\psi(\alpha)u = \alpha u$, where $\alpha u(t) = \alpha(t)u(t)$, for all $\alpha \in C_0(\R)$, $u \in L^2(\R)$ and $t \in \R$. Finally, let $\pi \colon \Z \rightarrow \calu(L^2(\R))$ be the map $(\pi(n)u) (t) = u(t-n)$ and consider the operator $\left(- i \frac{d}{dt}\right)$.
This operator is self-adjoint but \emph{not} bounded on $L^2(\R)$. We ``normalize'' it to obtain a bounded operator
$$
    T = \left(  \frac{x}{\sqrt{1+x^2}} \right) \left(- i \frac{d}{dt} \right)\,.
$$
This notation means that the function $\frac{x}{\sqrt{1+x^2}}$ is applied using functional calculus to the operator $\left(- i \frac{d}{dt} \right)$. Note that the operator $\left(- i \frac{d}{dt} \right)$ is essentially self adjoint. Thus the function $\frac{x}{\sqrt{1+x^2}}$ can be applied to the unique self-adjoint extension of $\left(- i \frac{d}{dt} \right)$.

Equivalently, $T$ can be constructed using Fourier transform. Let $\calm_x$ be the operator ``multiplication by $x$''
$$
    \calm_x(f(x)) = xf(x)\,.
$$
The Fourier transform $\calf$ converts $- i \frac{d}{dt}$ to $\calm_x$, i.e. there is a commutative diagram
$$
\xymatrix{
    L^2(\R) \ar[r]^-{\calf} \ar[d]_-{- i \frac{d}{dt}} & L^2(\R) \ar[d]^-{\calm_x}\\
    L^2(\R) \ar[r]_-{\calf} & L^2(\R)\,.
}
$$
Let $\calm_{\frac{x}{\sqrt{1+x^2}}}$ be the operator ``multiplication by $\frac{x}{\sqrt{1+x^2}}$''
$$
    \calm_{\frac{x}{\sqrt{1+x^2}}} (f(x)) =\frac{x}{\sqrt{1+x^2}} f(x)\,.
$$
$T$ is the unique bounded operator on $L^2(\R)$ such that the following diagram is commutative
$$
\xymatrix{
    L^2(\R) \ar[r]^-{\calf} \ar[d]_-{T} & L^2(\R) \ar[d]^-{\calm_{\frac{x}{\sqrt{1+x^2}}}}\\
    L^2(\R) \ar[r]_-{\calf} & L^2(\R)\,.
}
$$
Then we have an equivariant odd $K$-cycle $(L^2(\R), \psi, \pi, T) \in \cale^\Z_1(\R)$.
\end{example}
Let $X$ be a proper $G$-space with compact, second countable quotient $G\backslash X$ and $\cale_1^G(X)$ defined as above.
The \emph{equivariant $K$-homology group} $K_1^G(X)$ is defined as the quotient
$$
    K_1^G(X) = \cale_1^G(X) \big/ \sim,
$$
where $\sim$ represents \emph{homotopy}, in a sense that will be made precise later (Section \ref{section:HomotopyAndKKtheory}). It is an abelian group with addition and inverse given by
\begin{eqnarray*}
    (H,\psi,\pi,T) + (H',\psi',\pi',T') &=& (H\oplus H',\psi \oplus \psi',\pi \oplus \pi',T \oplus T'), \\
    -(H,\psi,\pi,T) &=& (H,\psi,\pi,-T).
\end{eqnarray*}
\begin{remark}
The $K$-cycles defined above differ slightly from the $K$-cycles used by Kasparov \cite{Kasparov88}. However, the abelian group $K_1^G(X)$ is isomorphic to the Kasparov group $KK_G^1(C_0(X),\C)$, where the isomorphism is given by the evident map which views one of our $K$-cycles as one of Kasparov's $K$-cycles. In other words, the $K$-cycles we are using are more special than the $K$-cycles used by Kasparov, however the obvious map of abelian groups is an isomorphism. 
\end{remark}

We define \emph{even} $K$-cycles in a similar way, just dropping the condition of $T$ being self-adjoint.
\begin{definition}
Let $X$ be a proper $G$-space with compact, second countable quotient space $G\backslash X$.
An \emph{equivariant even $K$-cycle} for $X$ is a 4-tuple $(H,\psi, \pi, T)$ such that:

\begin{itemize}
    \item $(H,\psi,\pi)$ is a representation of the $G$-\csa{} $C_0(X)$,
    \item $T \in \mathcal{L}(H)$,
    \item $\pi(g)T - T\pi(g) = 0$ for all $g \in G$,
    \item $\psi(\alpha)T-T\psi(\alpha) \in \mathcal{K}(H)$ for all $\alpha \in C_0(X)$,
    \item $\psi(\alpha)(I-T^*T)\in \mathcal{K}(H)$ for all $\alpha \in C_0(X)$,
    \item $\psi(\alpha)(I - TT^*) \in \mathcal{K}(H)$ for all $\alpha \in C_0(X)$.
\end{itemize}
\end{definition}
Write $\cale_0^G(X)$ for the set of such equivariant even $K$-cycles.
\begin{remark}\label{rmk:grading}
In the literature the definition is somewhat more complicated. In particular, the Hilbert space $H$ is required to be $\Z/2$-graded. However, at the level of abelian groups, the abelian group $K_0^G(X)$ obtained from the equivariant even $K$-cycles defined here will be isomorphic to the Kasparov group $KK^0_G(C_0(X),\C)$ \cite{Kasparov88}. 
More precisely, let $(H, \psi, \pi, T, \omega)$ be a $K$-cycle in Kasparov's sense, where $\omega$ is a $\Z/2$-grading of the Hilbert space $H=H_0 \oplus H_1$, $\psi = \psi_0 \oplus \psi_1$, $\pi = \pi_0 \oplus \pi_1$ and $T$ is self-adjoint but off-diagonal
$$
	T = \left( \begin{array}{cc} 0 & T_- \\ T_+ & 0 \end{array}\right).
$$
To define the isomorphism from $KK^0_G(C_0(X),\C)$ to $K_0^G(X)$, we map a Kasparov cycle $(H,\psi,\pi,T,\omega)$ to $(H',\psi',\pi',T')$ where
\begin{eqnarray*}
	H' &=& \ldots H_0 \oplus H_0 \oplus H_0 \oplus H_1 \oplus H_1 \oplus H_1 \ldots\\
	\psi' &=& \ldots \psi_0 \oplus \psi_0 \oplus \psi_0 \oplus \psi_1 \oplus \psi_1 \oplus \psi_1 \ldots\\
	\pi' &=&  \ldots \pi_0 \oplus \pi_0 \oplus \pi_0 \oplus \pi_1 \oplus \pi_1 \oplus \pi_1 \ldots\\
\end{eqnarray*}
and $T'$ is the obvious right-shift operator, where we use $T_+$ to map the last copy of $H_0$ to the first copy of $H_1$. The isomorphism from $\cale_0^G(X)$ to $KK^0_G(C_0(X),\C)$ is given by
$$
	(H, \psi, \pi, T) \mapsto (H\oplus H, \psi\oplus \psi, \pi\oplus \pi, \left( \begin{array}{cc} 0 & T^* \\ T & 0 \end{array}\right)).
$$
\end{remark}

Let $X$ be a proper $G$-space with compact, second countable quotient $G\backslash X$ and $\cale_0^G(X)$ as above.
The \emph{equivariant $K$-homology group} $K_0^G(X)$ is defined as the quotient
$$
    K_0^G(X) = \cale_0^G(X) \big/ \sim,
$$
where $\sim$ is \emph{homotopy}, in a sense that will be made precise later. It is an abelian group with addition and inverse given by
\begin{eqnarray*}
    (H,\psi,\pi,T) + (H',\psi',\pi',T') &=& (H\oplus H',\psi \oplus \psi',\pi \oplus \pi',T \oplus T'), \\
    -(H,\psi,\pi,T) &=& (H,\psi,\pi,T^*).
\end{eqnarray*}
\begin{remark}
Since the even $K$-cycles are more general, we have $\cale_1^G(X) \subset \cale_0^G(X)$. However, this inclusion induces the zero map from $K^G_1(X)$ to $K^G_0(X)$.
\end{remark}

\subsection{Functoriality}
Equivariant $K$-homology gives a (covariant) functor between the category proper $G$-spaces with compact quotient and the category of abelian groups. Indeed, given a continuous $G$-map $f \colon X \rightarrow Y$ between proper $G$-spaces with compact quotient, it induces a map $\widetilde{f} \colon C_0(Y) \rightarrow C_0(X)$ by $\widetilde{f}(\alpha) = \alpha \circ f$ for all $\alpha \in C_0(Y)$. Then, we obtain homomorphisms of abelian groups
$$
    K_j^G(X) \longrightarrow K_j^G(Y) \quad j = 0,1
$$
by defining, for each $(H,\psi,\pi,T) \in \cale^G_j(X)$,
$$
    (H,\psi,\pi,T) \mapsto (H,\psi \circ \widetilde{f},\pi,T)\,.
$$

\subsection{The index map}
Let $X$ be a proper second countable $G$-space with compact quotient $G\backslash X$. There is a map of abelian groups
\begin{eqnarray*}
  K_j^G(X) &\longrightarrow& K_j\left(C_r^*(G)\right) \nonumber \\
  (H,\psi, \pi, T) &\mapsto& \textup{Index}(T)
\end{eqnarray*}
for $j=0,1$. It is called the \emph{index map} and will be defined in Section \ref{section:IndexMap}.

This map is natural, that is, if $X$ and $Y$ are
proper second countable $G$-spaces with compact quotient and if $f \colon X
\rightarrow Y$ is a continuous $G$-equivariant map, then the
following diagram commutes:
$$ \xymatrix{
    K_j^G(X) \ar[r]^-{f_*} \ar[d]_-{\textup{Index}} & K_j^G(Y) \ar[d]^-{\textup{Index}}\\
    K_j^G(C_r^*(G)) \ar[r]^= & K_j^G(C_r^*(G)).
}$$
We would like to define equivariant $K$-homology and the index map for $\EG$.
However, the quotient of $\EG$ by the $G$-action might not be compact. The
solution will be to consider all proper second countable $G$-subspaces with compact
quotient.

\begin{definition}
Let $Z$ be a proper $G$-space. 
We call $\Delta \subseteq Z$ \emph{$G$-compact} if
\begin{enumerate}
  \item $gx \in \Delta$ for all $g \in G$, $x \in \Delta$,
  \item $\Delta$ is a proper $G$-space,
  \item the quotient space $G\backslash\Delta$ is compact.
\end{enumerate}
That is, $\Delta$ is a $G$-subspace which is proper as a $G$-space and has compact quotient $G\backslash \Delta$.
\end{definition}
\begin{remark}
Since we are always assuming that $G$ is locally compact, Hausdorff and second countable, we may also assume without loss of generality that any $G$-compact subset of $\EG$ is second countable. From now on we shall assume that $\EG$ has this property.
\end{remark}
We define the \emph{equivariant $K$-homology of \EG{}  with
$G$-compact supports} as the direct limit
$$
    K_j^G(\EG) = \lim_{\substack{\longrightarrow\\ \Delta \subseteq \EG\\ \textup{ $G$-compact}}} K_j^G(\Delta)\,.
$$
There is then a well-defined index map on the direct limit
\begin{eqnarray}\label{eqn:IndexMapEG}
  \mu \colon K_j^G(\EG) &\longrightarrow& K_j(C_r^* G) \nonumber \\
  (H,\psi,\pi,T) &\mapsto& \textup{Index}(T),
\end{eqnarray}
as follows.
Suppose that $\Delta \subset \Omega$ are $G$-compact. By the naturality of the functor $K^G_j(-)$, there is a commutative diagram
$$
\xymatrix{
    K^G_j(\Delta) \ar[r] \ar[d]_-{\textrm{Index}} & K^G_j(\Omega) \ar[d]^-{\textrm{Index}}\\
    K_j(C_r^* G) \ar[r]^= & K_j(C_r^* G)\,,
}
$$
and thus the index map is defined on the direct limit.





\section{The discrete case}\label{section:DiscreteCase}
We discuss several aspects of the Baum-Connes conjecture when the group is discrete.

\subsection{Equivariant $K$-homology}
For a discrete group $\Gamma$, there is a simple description of $K_j^\Gamma(\underline{E}\Gamma)$ up to torsion, in purely algebraic terms, given by a Chern character. Here we follow section 7 in \cite{BCH94}.

Let $\Gamma$ be a (countable) discrete group. Define $F\Gamma$ as the
set of finite formal sums
$$
    F\Gamma = \Bigg\{ \sum_{\textup{finite}} \lambda_\gamma [\gamma] \textup{ where } \gamma \in \Gamma, \textup{order}(\gamma) <
        \infty,\ \lambda_\gamma \in \C \Bigg\}\,.
$$
$F\Gamma$ is a complex vector space and also a $\Gamma$-module with
$\Gamma$-action:
$$
    g \cdot \left( \sum_{\lambda \in \Gamma} \lambda_\gamma [\gamma]
    \right) = \sum_{\lambda \in \Gamma} \lambda_\gamma [g\gamma
    g^{-1}]\,.
$$
Note that the identity element of the group has order 1 and therefore $F\Gamma \neq 0$.

Consider $H_j(\Gamma; F\Gamma)$, $j \ge 0$, the homology groups of
$\Gamma$ with coefficients in the $\Gamma$-module $F\Gamma$.

\begin{remark}
This is standard homological algebra, with no topology involved ($\Gamma$ is a discrete group and $F\Gamma$ is a
non-topologized module over $\Gamma$). They are classical homology
groups and have a purely algebraic description (cf.~\cite{Brown82}).
In general, if $M$ is a $\Gamma$-module then $H_*(\Gamma;M)$ is isomorphic to $H_*(B\Gamma;\underline{M})$, where $\underline{M}$ means the local system on $B\Gamma$ obtained from the $\Gamma$-module $M$.
\end{remark}

Let us write $K^{\textup{top}}_j(\Gamma)$ for
$K_j^\Gamma(\underline{E}\Gamma)$, $j=0,1$. There is a Chern
character $\textup{ch} \colon K^{\textup{top}}_*(\Gamma) \rightarrow
H_*(\Gamma; F\Gamma)$ which maps into odd, respectively even,
homology
$$
    \textup{ch} \colon K^{\textup{top}}_j(\Gamma) \rightarrow
    \bigoplus_{l\ge 0} H_{j+2l}(\Gamma; F\Gamma) \quad j=0,1.
$$
This map becomes an isomorphism when tensored with \C (cf.~\cite{BC88} or \cite{Luck02}).
\begin{proposition}\label{prop:TensoredChernCharacter}
The map
$$
    \textup{ch}\otimes_\Z \C \, \colon\; K^{\textup{top}}_j(\Gamma) \otimes_\Z \C \longrightarrow
    \bigoplus_{l \ge 0} H_{j+2l}(\Gamma; F\Gamma) \quad j=0,1
$$
is an isomorphism of vector spaces over \C.
\end{proposition}


\begin{remark}
If $G$ is finite, the rationalized Chern character becomes the
character map from $R(G)$, the complex representation ring of $G$,
to class functions, given by $\rho \mapsto \chi(\rho)$ in the even
case, and the zero map in the odd case.
\end{remark}


If the Baum-Connes conjecture is true for $\Gamma$, then Proposition
\ref{prop:TensoredChernCharacter} computes the tensored
topological $K$-theory of the reduced \csa{} of $\Gamma$.
\begin{corollary}
If the Baum-Connes conjecture is true for $\Gamma$ then
$$
    K_j(C_r^*\Gamma)\otimes_\Z \C \cong
    \bigoplus_{l \ge 0} H_{j+2l}(\Gamma; F\Gamma)\quad j=0,1\,.
$$
\end{corollary}
\subsection{Some results on discrete groups}
We recollect some results on discrete groups
which satisfy the Baum-Connes conjecture.
\begin{theorem}[N. Higson, G. Kasparov \cite{HK97}]
If $\Gamma$ is a discrete group which is \emph{amenable} (or, more
generally, \emph{a-T-menable}) then the Baum-Connes conjecture
is true for $\Gamma$.
\end{theorem}

\begin{theorem}[I.Mineyev, G. Yu \cite{MY02}; independently V. Lafforgue \cite{Lafforgue02}] If $\Gamma$ is a
discrete group which is \emph{hyperbolic} (in Gromov's sense) then
the Baum-Connes conjecture is true for $\Gamma$.
\end{theorem}

\begin{theorem}[Schick \cite{Schick07}] Let $B_n$ be the \emph{braid group} on $n$ strands, for any positive integer $n$. Then the Baum-Connes conjecture is true for $B_n$.
\end{theorem}

\begin{theorem}[Matthey, Oyono-Oyono, Pitsch \cite{MOP08}]
Let $M$ be a connected orientable 3-dimensional manifold
(possibly with boundary). Let $\Gamma$ be the fundamental group of
$M$. Then the Baum-Connes conjecture is true for $\Gamma$.
\end{theorem}

The Baum-Connes index map has been shown to be injective or rationally injective for some classes of groups. For example, it is injective for countable subgroups of $GL(n,K)$, $K$ any field \cite{GHW05}, and injective for
\begin{itemize}
\item closed subgroups of connected Lie groups \cite{Kasparov95};

\item closed subgroups of reductive $p$-adic groups \cite{KS91}.
\end{itemize}
More results on groups satisfying the Baum-Connes conjecture can be found in \cite{LR05}.

The Baum-Connes conjecture remains a widely open problem. For example, it is not known for $SL(n,\Z)$, $n\ge 3$. These infinite discrete groups have Kazhdan's property (T) and hence they are not a-T-menable. 
On the other hand, it is known that the index map is injective for $SL(n,\Z)$ (see above) and the groups $K^G_j(\EG)$ for $G=SL(3,\Z)$ have been calculated \cite{Sanchez-Garcia08}.

\begin{remark}
The conjecture might be too general to be true for all groups. Nevertheless, we expect it to be true for a large family of groups, in particular for all exact groups (a groups $G$ is \emph{exact} if the functor $C_r^*(G,-)$, as defined in \ref{section:ReducedCrossedProduct}, is exact).
\end{remark}

\subsection{Corollaries of the Baum-Connes Conjecture}
The Baum-Connes conjecture is related to a great number of
conjectures in functional analysis, algebra, geometry and topology.
Most of these conjectures follow from either the injectivity or the
surjectivity of the index map. A significant example is the Novikov
conjecture on the homotopy invariance of higher signatures of
closed, connected, oriented, smooth manifolds. This conjecture
follows from the injectivity of the rationalized index map
\cite{BCH94}. For more information on conjectures related to
Baum-Connes, see the appendix in \cite{Mislin03}.

\begin{remark}
By a ``corollary'' of the Baum-Connes conjecture we mean: if the Baum-Connes conjecture
is true for a group $G$ then the corollary is true for that group $G$. (For instance, in the Novikov conjecture $G$ is the fundamental group of the manifold.)
\end{remark}

\section{The compact case}\label{section:CompactCase}

If $G$ is compact, we can take \EG{} to be a one-point space. On the
other hand, $K_0 (C_r^*G) = R(G)$ the (complex) representation ring
of $G$, and $K_1 (C_r^*G) = 0$ (see Remark below). Recall that $R(G)$ is the
Grothendieck group of the category of finite dimensional (complex)
representations of $G$. It is a free abelian group with one
generator for each distinct (i.e.~non-equivalent) irreducible
representation of $G$. 
\begin{remark}
When $G$ is compact, the reduced \csa{} of $G$ is a direct sum (in the \csa{} sense) over the irreducible representations of $G$, of matrix algebras of dimension equal to the dimension of the representation. The $K$-theory functor commutes with direct sums and $K_j(M_n(\mathbb{C})) \cong K_j(\mathbb{C})$, which is $\mathbb{Z}$ for $j$ even and 0 otherwise (Theorem \ref{Thm:Bott2}).
\end{remark}

Hence the index map takes the form
$$
    \mu \colon K_G^0(point) \longrightarrow R(G)\,,
$$
for $j=0$ and is the zero map for $j=1$.

Given $(H,\psi,T,\pi) \in \cale^0_G(point)$, we may assume within the
equivalence relation on $\cale^0_G(point)$ that
\begin{gather*}
    \psi(\lambda) = \lambda I \quad \textup{ for all } \lambda \in C_0(point) = \C\,,
\end{gather*}
where $I$ is the identity operator of the Hilbert space $H$.
Hence the non-triviality
of $(H,\psi,T,\pi)$ is coming from
\begin{gather*}
    I - TT^* \in \calk(H)\,,\ \textup{ and }\
    I - T^*T \in \calk(H)\,,
\end{gather*}
that is, $T$ is a Fredholm operator.
Therefore
\begin{gather*}
    \dim_\C\left(\textup{ker}(T)\right) < \infty,\\
    \dim_\C\left(\textup{coker}(T)\right) < \infty,
\end{gather*}
hence $\textup{ker}(T)$ and $\textup{coker}(T)$ are finite
dimensional representations of $G$ (recall that $G$ is acting via
$\pi \colon G \rightarrow \call(H)$). Then
$$
    \mu(H,\psi,T,\pi) = \textup{Index}(T) = \textup{ker}(T) - \textup{coker}(T) \in
    R(G)\,.
$$

\begin{remark}
The assembly map for $G$ compact just described is an isomorphism (exercise).
\end{remark}

\begin{remark}
In general, for $G$ non-compact, the
elements of $K_0^G(X)$ can be viewed as generalized elliptic
operators on \EG, and the index map $\mu$ assigns to such an
operator its `index', $\textup{ker}(T) - \textup{coker}(T)$, in some
suitable sense \cite{BCH94}. This should be made precise later using Kasparov's descent map and an appropriate Kasparov product (Section \ref{section:IndexMap}).
\end{remark}


\section{Equivariant $K$-homology for \Gcsa{}s}\label{section:EquivariantKhomology2}
We have defined equivariant $K$-homology for $G$-spaces in Section \ref{section:EquivariantKHomology}. Now we define equivariant $K$-homology for a separable \Gcsa{} $A$ as the $KK$-theory groups $K^j_G(A,\C)$, $j=0,1$. This generalises the previous construction since $K_j^G(X) = KK_G^j(C_0(X), \mathbb{C})$. Later on we shall define $KK$-theory groups in full generality (Sections \ref{section:HomotopyAndKKtheory} and \ref{section:EquivariantKhomology3}).

\begin{definition} Let $A$ be a separable \Gcsa{}. Define $\cale^1_G(A)$ to be the set of 4-tuples
$$
    \{ (H,\psi,\pi,T) \}
$$
such that $(H,\psi,\pi)$ is a representation of the \Gcsa{} $A$, $T \in \call(H)$, and the following conditions are satisfied:
\begin{itemize}
  \item $T = T^*$,
  \item $\pi(g)T-T\pi(g)\in\calk(H)$,
  \item $\psi(a)T-T\psi(a)\in\calk(H)$,
  \item $\psi(a)(I-T^2)\in\calk(H)$,
\end{itemize}
for all $g \in G$, $a \in A$.
\end{definition}

\begin{remark}
Note that this is not quite $\cale^G_1(X)$ when $A=C_0(X)$ and $X$
is a proper $G$-space with compact quotient, since the third
condition is more general than before. However, the inclusion $\cale^G_1(X)
\subset \cale_G^1(C_0(X))$ gives an isomorphism of abelian groups so that $K_1^G(X) = KK_G^1(C_0(X), \C)$ (as defined below). The point is that, for a proper $G$-space with compact quotient, an averaging argument using a cut-off function and the Haar measure of the group $G$ allows us to assume that the operator $T$ is $G$-equivariant.
\end{remark}

Given a separable \Gcsa{} $A$, we define the $KK$-group $KK^1_G(A,\C)$ as $\cale^1_G(A)$ modulo an
equivalence relation called \emph{homotopy}, which will be made
precise later. Addition in $KK^1_G(A,\C)$ is given by direct sum
\[
(H,\psi,\pi,T) + (H',\psi',\pi',T') = (H\oplus H',\psi\oplus
\psi',\pi\oplus\pi', T\oplus T')
\]
and the negative of an element by
$$
    -(H,\psi,\pi,T) = (H,\psi,\pi,-T)\,.
$$
\begin{remark}
We shall later define $KK^1_G(A,B)$ for a separable \Gcsa{}s $A$ and an arbitrary \Gcsa{} $B$ (Section \ref{section:HomotopyAndKKtheory}).
\end{remark}

Let $A$, $B$ be separable \Gcsa{}s. A $G$-equivariant $\ast$-homomorphism
$\phi \colon A \rightarrow B$ gives a map $\cale^1_G(B)
\rightarrow \cale^1_G(A)$ by
$$
    (H,\psi,\pi,T) \mapsto (H,\psi \circ \phi,\pi,T)\,,
$$
and this induces a map $KK^1_G(B,\C) \rightarrow KK^1_G(A,\C)$.
That is, $KK^1_G(A,\C)$ is a contravariant functor in
$A$.

For the even case, the operator $T$ is not required to be self-adjoint. 

\begin{definition}
Let $A$ be a separable \Gcsa. Define $\cale^0_G(A)$ as the set of 4-tuples
$$
     \{ (H,\psi,\pi,T) \}
$$
such that $(H,\psi,\pi)$ is a representation of the \Gcsa{} $A$, $T \in
\call(H)$ and the following conditions are satisfied:
\begin{itemize}
  \item $\pi(g)T-T\pi(g)\in\calk(H)$,
  \item $\psi(a)T-T\psi(a)\in\calk(H)$,
  \item $\psi(a)(I-T^*T)\in\calk(H)$,
  \item $\psi(a)(I-TT^*)\in\calk(H)$,
\end{itemize}
for all $g \in G$, $a \in A$.
\end{definition}
\begin{remark}
Again, if $X$ is a proper $G$-space with compact quotient, the inclusion $\cale^G_0(X) \subset \cale^0_G(C_0(X))$ gives an isomorphism in $K$-homology, so we can write $K_0^G(X) = KK^0(C_0(X), \C)$ (as defined below). The issue of the $\Z/2$-grading (which is present in the Kasparov definition but not in our definition) is dealt with as in Remark \ref{rmk:grading}.
\end{remark}

We define the $KK$-groups $KK^0_G(A,\C)$ as $\cale^0_G(A)$ modulo an
equivalence relation called \emph{homotopy}, which will be made
precise later. Addition in $KK^1_G(A,\C)$ is given by direct sum
\[
(H,\psi,\pi,T) + (H',\psi',\pi',T') = (H\oplus H',\psi\oplus
\psi',\pi\oplus\pi',T\oplus T')
\]
and the negative of an element by
$$
    -(H,\psi,\pi,T) = (H,\psi,\pi,T^*)\,.
$$
\begin{remark}
We shall later define in general $KK^0_G(A,B)$ for a separable \Gcsa{}s $A$ and an arbitrary \Gcsa{} $B$ (Section \ref{section:EquivariantKhomology3}).
\end{remark}

Let $A$, $B$ be separable \Gcsa{}s. A $G$-equivariant $\ast$-homomorphism
$\phi \colon A \rightarrow B$ gives a map $\cale^0_G(B)
\rightarrow \cale^0_G(A)$ by
$$
    (H,\psi,\pi,T) \mapsto (H,\psi \circ \phi,\pi,T)\,,
$$
and this induces a map $KK^0_G(B,\C) \rightarrow KK^0_G(A,\C)$.
That is, $KK^0_G(A,\C)$ is a contravariant functor in $A$.

\section{The conjecture with coefficients}\label{section:BCCwC}
There is a generalized version of the Baum-Connes conjecture, known as the \emph{Baum-Connes conjecture with coefficients},
which adds coefficients in a \Gcsa{}. We recall the definition of \Gcsa.
\begin{definition}
A \Gcsa{} is a \csa{} $A$ with a given continuous action of $G$
$$G \times A \longrightarrow A$$
such that $G$ acts by $C^*$-algebra automorphisms. Continuity means that, for each $a \in A$, the map $G \rightarrow A$, $g \mapsto ga$ is a continuous map.
\end{definition}

\begin{remark}\label{rmk:trivialGcsa}
Observe that the only $\ast$-homomorphism of \C{} as a \csa{} is the
identity. Hence the only $G$-\csa{} structure on \C{} is the one
with trivial $G$-action.
\end{remark}

Let $A$ be a \Gcsa. Later we shall define the reduced crossed-product \csa{} $C_r^*(G,A)$, and the equivariant $K$-homology group with coefficients $K^G_j(\EG,A)$. These constructions reduce to $C_r^*(G)$, respectively $K^G_j(\EG)$, when $A=\C$. Moreover, the index map extends to this general setting and is also conjectured to be an isomorphism.
\begin{conjecture}[P.~Baum, A.~Connes, 1980]\label{conj:BCCwC} Let $G$ be a locally
compact, Hausdorff, second countable, topological group, and let $A$
be any \Gcsa{}. Then
$$
    \mu \colon K_j^G(\EG,A) \longrightarrow K_j (C_r^*(G,A)) \quad j=0,1
$$
is an isomorphism.
\end{conjecture}
Conjecture \ref{conj:BCC} follows as a particular case when $A=\C$. A fundamental difference is that the conjecture with coefficients is subgroup closed, that is, if it is true for a group $G$ for \emph{any} coefficients then it is true, for \emph{any} coefficients, for any closed subgroup of $G$.

The conjecture with coefficients has been proved for:
\begin{itemize}
   \item compact groups,
   \item abelian groups,
   \item groups acting simplicially on a tree with all vertex stabilizers satisfying the conjecture with coefficients \cite{Oyono-Oyono01},
   \item amenable groups and, more generally, a-T-menable groups (groups with the Haagerup property) \cite{HK01},
   \item the Lie group $Sp(n,1)$ \cite{Julg02},
   \item 3-manifold groups \cite{MOP08}.
\end{itemize}
For more examples of groups satisfying the conjecture with coefficients see \cite{LR05}.


\subsection*{Expander graphs}
Suppose that $\Gamma$ is a finitely generated, discrete group which
contains an expander family \cite{DSV03} in its Cayley graph as a
subgraph. Such a $\Gamma$ is a counter-example to the conjecture
with coefficients \cite{HLS02}. M.~Gromov outlined a proof that such
$\Gamma$ exists. A number of mathematicians are now filling in the
details. It seems quite likely that this group exists.

\section{Hilbert modules}\label{section:HilbertModules}
In this section we introduce the concept of Hilbert module over a \csa. It generalises the definition of Hilbert space by allowing the inner product to take values in a \csa. Our main application will be the definition of the reduced crossed-product \csa{} in Section \ref{section:ReducedCrossedProduct}. For a concise reference on Hilbert modules see \cite{Lance95}.

\subsection{Definitions and examples}
Let $A$ be a \csa.
\begin{definition}
An element $a \in A$ is \emph{positive} (notation: $a \ge 0$) if
there exists $b \in A$ with $bb^*=a$.
\end{definition}
The subset of positive elements, $A^+$, is a convex cone (closed under positive linear combinations) and $A^+ \cap (-A^+) = \{0\}$ \cite[1.6.1]{Dixmier77}. Hence we have a well-defined partial ordering in $A$ given by $x \ge y \iff x-y \ge 0$.
\begin{definition}
A \emph{pre-Hilbert $A$-module} is a right $A$-module $\calh$ with a
given $A$-valued inner product $\langle \ , \ \rangle$ such that:

\begin{itemize}
  \item $\langle u, v_1 + v_2 \rangle = \langle u, v_1 \rangle + \langle u, v_2 \rangle$,
  \item $\langle u, v a \rangle = \langle u, v \rangle a$,
  \item $\langle u, v \rangle = \langle v, u \rangle^*$,
  \item $\langle u, u \rangle \ge 0$,
  \item $\langle u, u \rangle = 0 \Leftrightarrow u=0$,
\end{itemize}
for all $u,v,v_1,v_2 \in \calh$, $a \in A$.
\end{definition}

\begin{definition}
A \emph{Hilbert $A$-module} is a pre-Hilbert $A$-module which is
complete with respect to the norm
$$
    \| u \| = \| \langle u, u \rangle \|^{1/2}\,.
$$
\end{definition}

\begin{remark}
If $\calh$ is a Hilbert $A$-module and $A$ has a unit $1_A$ then
$\calh$ is a complex vector space with
$$
    u \lambda = u (\lambda 1_A) \qquad u \in \calh, \lambda \in
    \C\,.
$$
If $A$ does not have a unit, then by using an
\emph{approximate identity} \cite{Pedersen79} in $A$, it is also a complex vector space.
\end{remark}

\begin{example}\label{eg:An}
Let $A$ be a \csa{} and $n \ge 1$. Then $A^n = A \oplus \ldots
\oplus A$ is a Hilbert $A$-module with operations
\begin{gather*}
    (a_1,\ldots,a_n) + (b_1,\ldots,b_n) = (a_1+b_1, \ldots,
    a_n+b_n),\\
    (a_1,\ldots,a_n) a = (a_1 a, \ldots,
    a_n a),\\
    \langle (a_1,\ldots,a_n), (b_1,\ldots,b_n) \rangle = a_1^*b_1 +
    \ldots + a_n^*b_n,\\
\end{gather*}
for all $a_j, b_j, a \in A$.
\end{example}

\begin{example}
Let $A$ be a \csa. Define
$$
    \calh = \bigg\{ (a_1, a_2, \ldots) \,\big|\, a_j \in A, \ \sum_{j=1}^\infty a_j^* a_j
    \textup{ is norm-convergent in $A$} \bigg\},
$$
with operations
\begin{gather*}
    (a_1,a_2,\ldots) + (b_1,b_2\ldots) = (a_1+b_1,a_2+b_2, \ldots),\\
    (a_1,a_2,\ldots) a = (a_1 a, a_2 a, \ldots),\\
    \langle (a_1,a_2,\ldots), (b_1,b_2\ldots) \rangle = \sum_{j=1}^\infty a_j^*b_j.
\end{gather*}
\end{example}

The previous examples can be generalized. Note that a \csa{} $A$ is a Hilbert module over itself with inner product $\langle a, b \rangle = a^*b$.

\begin{example}
If $\mathcal{H}_1, \ldots, \mathcal{H}_n$ are Hilbert $A$-modules then the direct sum $\mathcal{H}_1 \oplus \ldots \oplus \mathcal{H}_n$ is a Hilbert $A$-module with
$$
   \langle (x_1,\ldots, x_n), (y_1,\ldots, y_n) \rangle = \sum_i x_i^*y_i\,.
$$
We write  $\calh^n$ for the direct sum of $n$ copies of a Hilbert $A$-module $\calh$.
\end{example}

\begin{example}
If $\{\mathcal{H}_i\}_{i \in \mathbb{N}}$ is a countable family of Hilbert $A$-modules then
$$
    \calh = \bigg\{ (x_1, x_2, \ldots) \,\big|\, x_i \in \calh_i, \ \sum_{j=1}^\infty \langle x_j, x_j \rangle
    \textup{ is norm-convergent in $A$} \bigg\}
$$
is a Hilbert $A$-module with inner product $\langle x, y \rangle = \sum_{j=1}^\infty \langle x_j, y_j \rangle$.
\end{example}

The following is our key example.

\begin{example}
Let $G$ be a locally compact, Hausdorff, second countable,
topological group. Fix a left-invariant Haar measure $dg$ for $G$.
Let $A$ be a \Gcsa. Then $L^2(G,A)$ is a Hilbert $A$-module defined as follows. Denote by $C_c(G,A)$ the set of all continuous compactly supported functions from $G$ to $A$. On $C_c(G,A)$ consider the norm
$$
	\norm{f} = \left\|\int_G g^{-1}\left(f(g)^*f(g)\right)dg\right\|.
$$
$L^2(G,A)$ is the completion of $C_c(G,A)$ in this norm.
It is a Hilbert $A$-module with operations
\begin{gather*}
    (f+h)g = f(g) + h(g),\\
    (fa)g = f(g)(ga),\\
    \langle f, h \rangle = \int_G g^{-1} \left(f(g)^*h(g)\right)dg\,.
\end{gather*}
Note that when $A=\C$ the group action is trivial and we get $L^2(G)$ (cf. Remark \ref{rmk:trivialGcsa}).
\end{example}

\begin{definition}
Let $\calh$ be a Hilbert $A$-module. An $A$-module map $T \colon
\calh \rightarrow \calh$ is \emph{adjointable} if there exists an
$A$-module map $T^* \colon \calh \rightarrow \calh$ with
$$
    \langle Tu, v \rangle = \langle u, T^*v \rangle \quad \textup{for
    all } u, v \in \calh\,.
$$
\end{definition}
If $T^*$ exists, it is unique, and $\sup_{\|u\|=1} \|Tu\| < \infty$.
Set
$$
    \call(\calh) = \{ T \colon \calh \rightarrow \calh \,|\, T
    \textup{ is adjointable}\}\,.
$$
Then $\call(\calh)$ is a \csa{} with operations
\begin{gather*}
    (T+S)u=Tu+Su,\\
    (ST)u=S(Tu),\\
    (T\lambda)u=(Tu)\lambda\\
    T^* \textup{ as above,}\\
    \|T\| = \sup_{\|u\|=1} \|Tu\|.
\end{gather*}

\subsection{The reduced crossed-product
$C_r^*(G,A)$}\label{section:ReducedCrossedProduct}

Let $A$ be a
\Gcsa. Define
$$
    C_c(G,A) = \left\{ f \colon G \rightarrow A \,|\, f
        \textup{ continuous with compact support}\right\}\,.
$$
Then $C_c(G,A)$ is a complex algebra with
\begin{gather*}
    (f+h)g = f(g) + h(g),\\
    (f\lambda)g = f(g) \lambda,\\
    (f \ast h) (g_0) = \int_G f(g)\left(gh(g^{-1}g_0)\right)dg,
\end{gather*}
for all $g, g_0 \in G$, $\lambda \in \C$, $f, h \in C_c(G,A)$. The
product $\ast$ is called \emph{twisted convolution}.

Consider the Hilbert $A$-module $L^2(G,A)$. There is an injection of algebras
\begin{align*}
   C_c(G,A) &\hookrightarrow \call(L^2(G,A))\\
         f &\mapsto T_f
\end{align*}
where $T_f(u) = f \ast u$ is twisted convolution as above.
We define $C_r^*(G,A)$ as the \csa{} obtained by completing $C_c(G,A)$ with respect
to the norm $\|f\|=\|T_f\|$. When $A = \C$, the $G$-action must be trivial and $C_r^*(G,A) =
C_r^*(G)$.

\begin{example}
Let $G$ be a finite group, and $A$ a \Gcsa. Let $dg$ be the
Haar measure such that each $g \in G$ has mass 1. Then
$$
    C_r^*(G,A) = \bigg\{ \sum_{\gamma \in G} a_\gamma[\gamma]
    \,\Big|\;
    a_\gamma \in A \bigg\}
$$
with operations
\begin{gather*}
    \left( \sum_{\gamma \in G} a_\gamma[\gamma] \right) + \left(
    \sum_{\gamma \in G} b_\gamma[\gamma] \right) = \sum_{\gamma \in G} (a_\gamma +
    b_\gamma)[\gamma],\\
    \left( \sum_{\gamma \in G} a_\gamma[\gamma] \right) \lambda = \sum_{\gamma \in G} (a_\gamma
    \lambda)[\gamma],\\
    (a_\alpha[\alpha])(b_\beta[\beta]) = a_\alpha(\alpha
    b_\beta)[\alpha \beta] \quad \textup{(twisted convolution)},\\
    \left( \sum_{\gamma \in G} a_\gamma[\gamma] \right)^* = \sum_{\gamma \in G} (\gamma^{-1}
    a_\gamma^*)[\gamma^{-1}].
\end{gather*}
\end{example}
Here $a_\gamma[\gamma]$ denotes the function from $G$ to $A$ which has the value $a_\gamma$ at $\gamma$ and 0 at $g \neq \gamma$.

Let $X$ be a Hausdorff, locally compact $G$-space. We know that
$C_0(X)$ becomes a \Gcsa{} with $G$-action
$$
    (gf)(x) = f(g^{-1}x),
$$
for $g \in G$, $f \in C_0(X)$ and $x \in X$.
The reduced crossed-product $C^*_r(G,C_0(X))$ will be denoted
$C^*_r(G,X)$.

A natural question is to calculate the $K$-theory of this \csa{}. 
If $G$ is compact, this is the Atiyah-Segal group
$K_G^j(X)$, $j=0,1$. Hence for $G$ non-compact, $K_j(C_r^*(G,X))$ is
the natural extension of the Atiyah-Segal theory to the case when
$G$ is non-compact.


\begin{definition}
We call a $G$-space \emph{$G$-compact}
if the quotient space $G\backslash X$ (with the quotient topology) is compact.
\end{definition}
Let $X$ be a proper, $G$-compact $G$-space. Then a $G$-equivariant
$\C$-vector bundle $E$ on $X$ determines an element
$$
    [E] \in K_0 (C^*_r(G,X))\,.
$$
\begin{remark}
From $E$, a Hilbert module over $C^*_r(G,X)$ is constructed. This Hilbert $C^*_r(G,X)$-module determines an element in $KK_0(\mathbb{C}, C_r^*(G,X)) \cong K_0(C_r^*(G,X))$. Note that, quite generally, a Hilbert $A$-module determines an element in $KK_0(A)$ if and only if it is finitely generated.
\end{remark}

Recall that a $G$-equivariant vector bundle $E$ over $X$ is a
(complex) vector bundle $\pi \colon E \rightarrow X$ together with a
$G$-action on $E$ such that $\pi$ is $G$-equivariant and, for each
$p \in X$, the map  on the fibers $E_p \rightarrow E_{gp}$ induced
by multiplication by $g$ is linear.

\begin{theorem}[W.~L\"uck and B.~Oliver \cite{LO01}] If $\Gamma$ is a
(countable) discrete group and $X$ is a proper $\Gamma$-compact
$\Gamma$-space, then
$$
    K_0 (C_r^*(\Gamma, X)) = \textup{Grothendieck group of
    $\Gamma$-equivariant \C-vector bundles on $X$.}
$$
\end{theorem}
\begin{remark}
In \cite{LO01} this theorem is not explicitly stated. However, it follows from their results. For clarification see \cite{BHS} or \cite{EM09}.
\end{remark}

\begin{remark}
Let $X$ be a proper $G$-compact $G$-space. Let $\mathbb{I}$ be the
trivial $G$-equivariant complex vector bundle on $X$,
$$
    \mathbb{I} = X \times \C\,, \quad g(x,\lambda) = (gx, \lambda),
$$
for all $g \in G$, $x \in X$ and $\lambda \in \C$.
Then $[\,\mathbb{I}\,] \in K_0 (C_r^*(G,X))$.
\end{remark}

\subsection{Push-forward of Hilbert modules}
Let $A$, $B$ be \csa{}s, $\varphi
\colon A \rightarrow B$ a $\ast$-homomorphism and $\calh$ a Hilbert
$A$-module. We shall define a Hilbert $B$-module $\calh \otimes_A
B$, called the \emph{push-forward of $\calh$ with respect to $\varphi$} or \emph{interior tensor product} (\cite[Chapter 4]{Lance95}). First, form the algebraic tensor product $\calh \odot_A B =
\calh \otimes_A^{alg} B$ ($B$ is an $A$-module via $\varphi$). This
is an abelian group and also a (right) $B$-module
$$
    (h \otimes b)b' = h \otimes bb' \quad \textup{for all } h \in
    \calh, b,b' \in B\,.
$$
Define a $B$-valued inner product on $\calh \odot_A B$ by
$$
    \langle h \otimes b, h' \otimes b' \rangle = b^* \varphi(\langle h,
    h' \rangle) b'.
$$
Set
$$
    \caln = \{ \xi \in \calh \odot_A B \,|\, \langle \xi, \xi
    \rangle = 0 \}\,.
$$
$\caln$ is a $B$-sub-module of $\calh \odot_A B$ and $(\calh \odot_A B)/\caln$ is a
pre-Hilbert $B$-module.

\begin{definition}
Define $\calh \otimes_A B$ to be the Hilbert $B$-module obtained by
completing $(\calh \odot_A B)/\caln$.
\end{definition}

\begin{example} Let $X$ be a locally compact, Hausdorff space. Let $A =
C_0(X)$, $B=\C$ and $ev_p \colon C_0(X) \rightarrow \C$ the
evaluation map at a point $p \in X$. Then we can consider the
push-forward of a Hilbert $C_0(X)$-module $\calh$. This gives a
Hilbert space $\calh_p$. These Hilbert spaces do not form a vector bundle but something more general (not necessarily locally trivial),
sometimes called \emph{continuous field of Hilbert spaces} \cite[chapter 10]{Dixmier77}.
\end{example}

\section{Homotopy made precise and $KK$-theory}\label{section:HomotopyAndKKtheory}
We first define homotopy and Kasparov's $KK$-theory in the non-equivariant setting, for pairs of separable \csa{}s. 
A first introduction to $KK$-theory and further references can be found in \cite{Higson90}.  

Let $A$ be a \csa{} and let $\calh$ be a Hilbert $A$-module.
Consider $\call(\calh)$ the bounded operators on $\calh$. For each
$u,v \in \calh$ we have a bounded operator $\theta_{u,v}$ defined as
$$
    \theta_{u,v}(\xi) = u \langle v, \xi \rangle\,.
$$
It is clear that $\theta_{u,v}^* = \theta_{v,u}$. The $\theta_{u,v}$
are called \emph{rank one operators} on $\calh$. A \emph{finite
rank operator} on $\calh$ is any $T \in \call(\calh)$ such that $T$
is a finite sum of rank one operators,
$$
    T = \theta_{u_1,v_1} + \ldots + \theta_{u_n,v_n}\,.
$$
Let $\calk(\calh)$ be the closure (in $\call(\calh)$) of the set of
finite rank operators. $\calk(\calh)$ is an ideal in $\call(\calh)$.
When $A=\C$, $\calh$ is a Hilbert space and $\calk(\calh)$ coincides
with the usual compact operators on $\calh$.

\begin{definition}
$\calh$ is \emph{countably generated} if in $\calh$ there is a
countable (or finite) set such that the $A$-module generated by this
set is dense in $\calh$.
\end{definition}

\begin{definition}
Let $\calh_0$, $\calh_1$ be two Hilbert $A$-modules. We say that
$\calh_0$ and $\calh_1$ are \emph{isomorphic} if there exists an
$A$-module isomorphism $\Phi \colon \calh_0 \rightarrow \calh_1$
with
$$
    \langle u, v \rangle_0 = \langle \Phi u, \Phi v \rangle_1 \quad
    \textup{for all } u, v \in \calh_0\,.
$$
\end{definition}

We want to define non-equivariant $KK$-theory for pairs
of \csa{}s. Let $A$ and $B$ be \csa{}s where $A$ is also separable. Define the set
$$
    \cale^1(A,B) = \{ (\calh, \psi, T) \}
$$
such that $\calh$ is a countably generated Hilbert $B$-module, $\psi
\colon A \rightarrow \call(\calh)$ is a $\ast$-homomorphism, $T \in
\call(\calh)$, and the following conditions are satisfied:
\begin{itemize}
  \item $T = T^*$,
  \item $\psi(a)T - T\psi(a) \in \calk(\calh)$,
  \item $\psi(a) (I - T^2) \in \calk(\calh)$,
\end{itemize}
for all $a \in A$. We call such triples \emph{odd bivariant
$K$-cycles}.

\begin{remark}
In the Kasparov definition of $KK^1(A,B)$ \cite{Kasparov80}, the conditions of the $K$-cycles are the same as our conditions except that the requirement $T=T^*$ is replaced by $\psi(a)(T-T^*) \in \calk(H)$ for all $a \in A$. The isomorphism of abelian groups from the group defined using these bivariant $K$-cycles to the group defined using our bivariant $K$-cycles is obtained by sending a Kasparov cycle $(H,\psi,T)$ to $(H, \psi, \frac{T+T^*}{2})$.
\end{remark}

We say that two such triples $(\calh_0, \psi_0, T_0)$ and $(\calh_1,
\psi_1, T_1)$ in $\cale^1(A,B)$ are \emph{isomorphic} if there is an
isomorphism of Hilbert $B$-modules $\Phi \colon \calh_0 \rightarrow
\calh_1$ with
\begin{gather*}
    \Phi \psi_0(a) = \psi_1(a) \Phi,\\
    \Phi T_0 = T_1 \Phi,
\end{gather*}
for all $a \in A$. That is, the following diagrams commute
$$ \xymatrix{
    \calh_0 \ar[r]^-{\psi_0(a)} \ar[d]_-{\Phi} & \calh_0 \ar[d]^-{\Phi}\\
    \calh_1 \ar[r]_-{\psi_1(a)} & \calh_1
}\qquad \xymatrix{
    \calh_0 \ar[r]^-{T_0} \ar[d]_-{\Phi} & \calh_0 \ar[d]^-{\Phi}\\
    \calh_1 \ar[r]_-{T_1} & \calh_1
}$$

Let $A$, $B$, $D$ be \csa{}s where $A$ is also separable. A $\ast$-homomorphism
$\varphi \colon B \rightarrow D$ induces a map $\varphi_* \colon
\cale^1(A,B) \rightarrow \cale^1(A,D)$ by
$$
    \varphi_*(\calh,\psi,T) = (\calh \otimes_B D, \psi \otimes_B I,
    T \otimes_B I)
$$
where $I$ is the identity operator on $D$, that is,
$I(\alpha)=\alpha$ for all $\alpha \in D$.

We can now make the definition of homotopy precise. Consider the
\csa{} of continuous functions $C([0,1],B)$, and set $\rho_0$,
$\rho_1$ to be the $*$-homomorphisms
$$
\xymatrix@R=1ex{ C([0,1],B) \ar@<1ex>[r]^-{\rho_0}
\ar@<-1ex>[r]_-{\rho_1} & B}
$$
defined by $\rho_0(f)=f(0)$ and $\rho_1(f)=f(1)$. In particular, we
have induced maps
$$
    (\rho_j)_* \colon \cale^1(A,C([0,1],B)) \longrightarrow
    \cale^1(A,B) \qquad j=0,1
$$
for any separable \csa{} $A$.
\begin{definition}
Two triples $(\calh_0,\psi_0,T_0)$ and $(\calh_1,\psi_1,T_1)$ in
$\cale^1(A,B)$ are \emph{homotopic} if there exists $(\calh,\psi,T)$
in $\cale^1(A, C([0,1],B))$ with
$$
    (\rho_j)_*(\calh,\psi,T) \cong (\calh_j,\psi_j,T_j) \quad
    j=0,1\,.
$$
\end{definition}

The even case is analogous, removing the self-adjoint condition
$T=T^*$. 
\begin{remark}
As above, we do not require the Hilbert $B$-module $\calh$ to be $\Z/2$-graded. The isomorphism between the abelian group we are defining and the group $KK_0(A,B)$ as defined by Kasparov \cite{Kasparov80} is dealt with as before (see Remark \ref{rmk:grading}).
\end{remark}

Hence we have the set of \emph{even bivariant $K$-cycles}
$$
    \cale^0(A,B) = \{ (\calh, \psi, T) \}
$$
where $\calh$ is a countably generated Hilbert $B$-module, $\psi
\colon A \rightarrow \call(\calh)$ a $\ast$-homomorphism, $T \in
\call(\calh)$, and the following conditions are satisfied:
\begin{itemize}
  \item $\psi(a)T - T\psi(a) \in \calk(\calh)$,
  \item $\psi(a) (I - T^*T) \in \calk(\calh)$,
  \item $\psi(a) (I - TT^*) \in \calk(\calh)$,
\end{itemize}
for all $a \in A$. The remaining definitions carry over, in
particular the definition of \emph{homotopy} in $\cale^0(A,B)$.

We define the \emph{(non-equivariant) Kasparov $KK$-theory groups}
of the pair $(A,B)$ as
\begin{eqnarray*}
  KK^1(A,B) &=& \cale^1(A,B) / \textup{(homotopy)}, \\
  KK^0(A,B) &=& \cale^0(A,B) / \textup{(homotopy)}.
\end{eqnarray*}
A key property is that $KK$-theory incorporates $K$-theory of \csa{}s: for any \csa{} $B$, $KK^j(\mathbb{C},B)$ is isomorphic to $K_j(B)$ (see Theorem 25 in \cite{Meyer}).

\section{Equivariant $KK$-theory}\label{section:EquivariantKhomology3}
We generalize $KK$-theory to the equivariant setting. An alternative definition to ours, by means of a universal property, is described in Section 2 of Meyer's notes \cite{Meyer}. 

All through this section, let $A$ be a \Gcsa.
\begin{definition}
A $G$-Hilbert $A$-module is a Hilbert $A$-module $\calh$ with a
given continuous action
\begin{eqnarray*}
  G \times \calh &\rightarrow& \calh \\
  (g,v) &\mapsto& gv
\end{eqnarray*}
such that
\begin{enumerate}
  \item $g(u+v) = gu + gv$,
  \item $g(ua) = (gu)(ga)$,
  \item $\langle gu, gv \rangle = g \langle u, v \rangle$,
\end{enumerate}
for all $g \in G$, $u,v \in \calh$, $a \in A$.
\end{definition}
Here `continuous' means that for each $u \in \calh$, the map $G
\rightarrow \calh$, $g \mapsto gu$ is continuous.

\begin{example} If $A=\C$, a $G$-Hilbert $\C$-module is just a unitary
representation of $G$ (the action of $G$ on \C{} must be trivial).
\end{example}

\begin{remark}
Let $\calh$ be a $G$-Hilbert $A$-module. For
each $g \in G$, denote by $L_g$ the map
$$
    L_g \colon \calh \rightarrow \calh, \quad L_g(v) = gv\,.
$$
Note that $L_g$ might not be in $\call(\calh)$. But if $T \in
\call(\calh)$, then $L_g T L_g^{-1} \in \call(\calh)$.

%
Hence $G$ acts on the \csa{} $\call(\calh)$ by
$$
    g T = L_g T L_g^{-1}\,.
$$
\end{remark}

\begin{example}
Let $A$ be a \Gcsa. Set $n \ge 1$. Then $A^n$ is a $G$-Hilbert
$A$-module (cf.~Example \ref{eg:An}) with
$$
    g(a_1, \ldots, a_n) = (ga_1,\ldots,ga_n).
$$
\end{example}


Let $A$ and $B$ be \Gcsa{}s, where $A$ is also separable.
Define the set
$$
    \cale^0_G(A,B) = \{ (\calh, \psi, T) \}
$$
such that $\calh$ is a countably generated $G$-Hilbert $B$-module,
$\psi \colon A \rightarrow \call(\calh)$ is a $\ast$-homomorphism
with
$$
    \psi(ga) = g \psi(a) \quad \textup{for all } g \in G, a \in A\,,
$$
and $T \in \call(\calh)$, and so that the following conditions are
satisfied:
\begin{itemize}
  \item $gT - T \in \calk(\calh)$,
  \item $\psi(a)T - T\psi(a) \in \calk(\calh)$,
  \item $\psi(a) (I - T^*T) \in \calk(\calh)$,
  \item $\psi(a) (I - TT^*) \in \calk(\calh)$,
\end{itemize}
for all $g \in G$, $a \in A$.
We define
$$
    KK_G^0(A,B) = \cale^0_G(A,B) / \textup{(homotopy)}\,.
$$
The definition of \emph{homotopy} in Section
\ref{section:HomotopyAndKKtheory} can be defined in a straightforward way in this setting.

$KK^0_G(A,B)$ is an abelian group with addition and
negative
\begin{eqnarray*}
  (\calh,\psi,T) + (\calh',\psi',T') &=& (\calh\oplus \calh',\psi\oplus\psi',T\oplus T'), \\
  -(\calh,\psi,T) &=& (\calh,\psi,T^*)\,.
\end{eqnarray*}

The odd case is similar, just restricting to self-adjoint operators.
Define the set
$$
    \cale^1_G(A,B) = \{ (\calh, \psi, T) \}
$$
such that $\calh$ is a countably generated $G$-Hilbert $B$-module,
$\psi \colon A \rightarrow \call(\calh)$ is a $\ast$-homomorphism
with
$$
    \psi(ga) = g \psi(a) \quad \textup{for all } g \in G, a \in A\,,
$$
and $T \in \call(\calh)$, and so that the following conditions are
satisfied:
\begin{itemize}
  \item $T = T^*$,
  \item $gT - T \in \calk(\calh)$,
  \item $\psi(a)T - T\psi(a) \in \calk(\calh)$,
  \item $\psi(a) (I - T^2) \in \calk(\calh)$,
\end{itemize}
for all $g \in G$, $a \in A$. 

We define
$$
    KK_G^1(A,B) = \cale^1_G(A,B) / \textup{(homotopy)}\,.
$$
$KK^1_G(A,B)$ is an abelian group with addition and inverse given by
\begin{eqnarray*}
  (\calh,\psi,T) + (\calh',\psi',T') &=& (\calh\oplus \calh',\psi\oplus\psi',T\oplus T'), \\
  -(\calh,\psi,T) &=& (\calh,\psi,-T).
\end{eqnarray*}

\begin{remark}
In the even case we are not requiring a $\Z/2$-grading. The isomorphism to the abelian group defined by Kasparov \cite{Kasparov88} is given as in Remark \ref{rmk:grading}. Our general principle is that the even and odd cases are identical except that in the odd case the operator $T$ is required to be self-adjoint but not in the even case.
\end{remark}

Using equivariant $KK$-theory, we can introduce coefficients for equivariant $K$-homology. Let $X$ be a proper $G$-space with compact quotient. Recall that 
\begin{eqnarray*}
   K^G_j(X) &=& KK^j_G(C_0(X),\C) \quad \textup{ and}\\
   K_j^G(\EG) &=& \lim_{\substack{\longrightarrow\\ \Delta \subseteq \EG\\ \textup{ $G$-compact}}} K_j^G(\Delta)\,.
\end{eqnarray*}
We define the \emph{equivariant $K$-homology of $X$}, respectively \emph{of \EG}, with coefficients in a \Gcsa{} $A$ as
\begin{eqnarray*}
   K^G_j(X,A) &=& KK^j_G(C_0(X),A),\\
   K_j^G(\EG,A) &=& \lim_{\substack{\longrightarrow\\ \Delta \subseteq \EG\\ \textup{ $G$-compact}}} K_j^G(\Delta,A)\,.
\end{eqnarray*}




\section{The index map}\label{section:IndexMap}
Our definition of the index map uses the Kasparov product and the descent map.

\subsection{The Kasparov product}
Let $A$, $B$, $D$ be (separable) \Gcsa{}s. There is a product
$$
    KK_G^i(A,B) \otimes_\Z KK_G^j(B,D) \longrightarrow
    KK_G^{i+j}(A,D)\,.
$$
The definition is highly non-trivial. Some motivation and examples, in the non-equivariant case, can be found in \cite[Section 5]{Higson90}.

\begin{remark}
Equivariant $KK$-theory can be regarded as a category with objects
separable \Gcsa{}s and morphisms $\textup{mor}(A,B) = KK^i_G(A,B)$
(as a $\Z/2$-graded abelian group), and composition given by the
Kasparov product (cf.~\cite[Thm.~33]{Meyer}).
\end{remark}

\subsection{The Kasparov descent map}
Let $A$ and $B$ be (separable) \Gcsa{}s. There is a map between
the equivariant $KK$-theory of $(A,B)$ and the non-equivariant
$KK$-theory of the corresponding reduced crossed-product \csa{}s,
$$
    KK_G^j(A,B) \longrightarrow KK^j\left(C^*_r(G,A),
    C^*_r(G,B)\right) \quad j=0,1\,.
$$
The definition is also highly non-trivial and can be found in \cite[Section 3]{Kasparov88}. Alternatively, see Proposition 26 in Meyer's notes \cite{Meyer}.


\subsection{Definition of the index map}
We would like to define the index map
$$
    \mu \colon K_j^G(\EG) \longrightarrow K_j (C_r^*G)\,.
$$
Let $X$ be a proper $G$-compact $G$-space. First, we define a map
$$
    \mu \colon K_j^G(X) = KK_G^j(C_0(X),\C) \longrightarrow K_j (C_r^*G)
$$
to be the composition of the Kasparov descent map
    $$
        KK_G^j(C_0(X),\C) \longrightarrow KK^j\left(C^*_r(G,X),
        C^*_r(G)\right)\,
    $$
(the trivial action of $G$ on \C{} gives the crossed-product
$C_r^*(G,\C) = C_r^*G$) and the Kasparov product with the trivial
bundle
  $$
    \mathbb{I} \in K_0 (C_r^*(G,X)) = KK^0(\C,C_r^*(G,X)),
  $$
that is, the Kasparov product with the trivial vector bundle
$\mathbb{I}$, when $A=\C$, $B=C_r^*(G,X)$,
  $D=C_r^*G$ and $i=0$.

Recall that
$$
    K_j^G(\EG) = \lim_{\substack{\longrightarrow\\ \Delta \subset \EG\\ \textup{ $G$-compact}}} KK^j_G\left(C_0(\Delta),\C\right)\,.
$$

For each $G$-compact $\Delta \subset \EG$, we have a map as before
$$
    \mu \colon KK^j_G(C_0(\Delta),\C) \longrightarrow K_j (C_r^* G)\,.
$$
If $\Delta$ and $\Omega$ are two $G$-compact subsets of \EG{} with
$\Delta \subset \Omega$, then by naturality the following diagram
commutes:
$$ \xymatrix{
    KK^j_G(C_0(\Delta),\C) \ar[r] \ar[d] & KK^j_G(C_0(\Omega),\C) \ar[d]\\
    K_j C_r^*G \ar[r]^-= & K_j C_r^*G .
}
$$
Thus we obtain a well-defined map on the direct limit $\mu \colon
K_j^G(\EG) \rightarrow K_j C_r^*G$.

\subsection{The index map with coefficients}
The coefficients can be introduced in $KK$-theory at once.
Let $A$ be a \Gcsa. We would like to define the index map
$$
    \mu \colon K_j^G(\EG;A) \longrightarrow K_j C_r^*(G,A)\,.
$$
Let $X$ be a proper $G$-compact $G$-space and $A$ a \Gcsa. First, we define a map
$$
    \mu \colon KK_G^j(C_0(X),A) \longrightarrow K_j C_r^*(G,A)
$$
to be the composition of the Kasparov descent map
    $$
        KK_G^j(C_0(X),A) \longrightarrow KK^j\left(C^*_r(G,X),
        C^*_r(G,A)\right)\,
    $$
and the Kasparov product with the trivial bundle
  $$
    \mathbb{I} \in K_0 C_r^*(G,X) = KK^0(\C,C_r^*(G,X)).
  $$
For each $G$-compact $\Delta \subset \EG$, we have a map as above
$$
    \mu \colon KK^j_G(C_0(\Delta),A) \longrightarrow K_j C_r^* (G,A)\,.
$$
If $\Delta$ and $\Omega$ are two $G$-compact subsets of \EG{} with
$\Delta \subset \Omega$, then by naturality the following diagram commutes:
$$ \xymatrix{
    KK^j_G(C_0(\Delta),A) \ar[r] \ar[d] & KK^j_G(C_0(\Omega),A) \ar[d]\\
    K_j C_r^*(G,A) \ar[r]^-= & K_j C_r^*(G,A).
}
$$
Thus we obtain a well-defined map on the direct limit $\mu \colon K_j^G(\EG;A) \rightarrow K_j
C_r^*(G,A)$.

\section{A brief history of $K$-theory}\label{section:History}

\subsection{The $K$-theory genealogy tree}
Grothendieck invented $K$-theory to give a conceptual proof of the Hirzebruch--Riemann--Roch theorem. The subject has since then evolved in different directions, as summarized by the following diagram.

\vspace{1cm}

\xymatrix{\textup{A.~Grothendieck} &
\framebox{\textbf{Riemann-Roch}}\ar[d]^-{\textup{Atiyah +
Hirzebruch}}
& \txt{$K$-theory in\\ algebraic geometry}\\
\textup{J.~F.~Adams} & \framebox{\txt{{\bf Vector fields}\\ {\bf on spheres}}}
\ar[dl]^-{\txt{{\phantom{x} A.~Connes}}} \ar[dr]_-{\txt{{H.~Bass}\\
{D.~Quillen}\\ {J.~H.~C.~Whitehead}}}
& \txt{$K$-theory\\ in topology}\\
\framebox{\txt{{\bf $K$-theory for}\\ {\bf \csa{}s}}} & &
\framebox{\txt{{\bf Algebraic}\\ {\bf $K$-theory}}}}

\vspace{0.6cm}
Atiyah and Hirzebruch defined topological $K$-theory. J.~F.~Adams then used the Atiyah-Hirzebruch theory to solve the problem of vector fields on spheres. \csa{} $K$-theory developed quite directly out of Atiyah-Hirzebruch topological $K$-theory. From its inception, \csa{} $K$-theory has been closely linked to problems in geometry-topology (Novikov conjecture, Gromov-Lawson-Rosenberg conjecture, Atiyah-Singer index theorem) and to classification problems within \csa{}s. More recently, \csa{} $K$-theory has played an essential role in the new subject of non-commutative geometry.

Algebraic $K$-theory was a little slower to develop \cite{Weibel99}; much of the early development in the 1960s was due to H.~Bass, who organized the theory on $K_0$ and $K_1$ and defined the negative $K$-groups. J.~Milnor introduced $K_2$. Formulating an appropriate definition for higher algebraic $K$-theory proved to be a difficult and elusive problem. Definitions were proposed by several authors, including J.~Milnor and Karoubi-Villamayor. A remarkable breakthrough was achieved by D.~Quillen with his plus-construction. The resulting definition of higher algebraic $K$-theory (i.e.~Quillen's algebraic $K$-theory) is perhaps the most widely accepted today. Many significant problems and results (e.g.~the Lichtenbaum conjecture) have been stated within the context of Quillen algebraic $K$-theory. In some situations, however, a different definition is relevant. For example, in the recently proved Bloch-Kato conjecture, it is J.~Milnor's definition of higher algebraic K-theory which is used.

Since the 1970s, $K$-theory has grown considerably, and its connections with other parts of mathematics have expanded. For the interested reader, we have included a number of current $K$-theory textbooks in our reference list (\cite{Blackadar98}, \cite{RLL00}, \cite{Rosenberg94}, \cite{Srinivas96}, \cite{Wegge-Olsen93}, \cite{WeibelK-Book}). For a taste of the current developments, it is useful to take a look at the \emph{Handbook of $K$-theory} \cite{FG05} or at the lectures in \cite{sedanoproceedings}. The \emph{Journal of $K$-theory} (as well as its predecessor, \emph{$K$-theory}) is dedicated to the subject, as is the website maintained by D.~Grayson at \verb+http://www.math.uiuc.edu/K-theory+. This site, started in 1993, includes a preprint archive which at the moment when this is being written contains 922 preprints. Additionally, see the \emph{Journal of Non-Commutative Geometry} for current results involving \csa{} $K$-theory.

Finally, we have not in these notes emphasized cyclic homology. However, cyclic (co-)homology is an allied theory to $K$-theory and any state-of-the-art survey of $K$-theory would have to recognize this central fact.

\subsection{The Hirzebruch--Riemann--Roch theorem}
Let $M$ be a non-singular projective algebraic variety over \C. Let
$E$ be an algebraic vector bundle on $M$. Write $\underline{E}$ for
the sheaf (of germs) of algebraic sections of $E$. For each $j \ge
0$, consider $H^j(M,\underline{E})$ the $j$-th cohomology group of $M$
using $\underline{E}$.

\begin{lemma}
For all $j\ge 0$, $\dim_\C H^j(M, \underline{E}) < \infty$ and for
$j > \dim_\C(M)$, $H^j(M, \underline{E}) = 0$.
\end{lemma}
Define the \emph{Euler characteristic} of $M$ with respect to $E$ as
$$
    \chi(M,E) = \sum_{j=0}^n (-1)^j \dim_\C H^j(M, \underline{E})\,,
    \quad \textup{where } n=\dim_\C (M)\,.
$$

\begin{theorem}[Hirzebruch--Riemann--Roch] Let $M$ be a non-singular
projective algebraic variety over \C{} and let $E$ be an algebraic
vector bundle on $M$. Then
$$
    \chi(M,E) = (\textup{ch}(E) \cup \textup{Td}(M))[M]
$$
where $\textup{ch}(E)$ is the Chern character of $E$, $\textup{Td}(M)$ is the Todd class of $M$ and $\cup$ stands for the cup product.
\end{theorem}











\subsection{The unity of $K$-theory}\label{section:Comparison}
We explain how $K$-theory for \csa{}s is a particular case of
algebraic $K$-theory of rings.

Let $A$ be a \csa. Consider the inclusion
\begin{eqnarray}\label{eqn:Inclusion}
  M_n(A) &\hookrightarrow & M_{n+1}(A) \nonumber\\
  \left( \begin{array}{ccc}
            a_{11} & \ldots & a_{1n}\\
            \vdots &        & \vdots\\
            a_{n1} & \ldots & a_{nn}\\
  \end{array} \right)
   & \mapsto &
   \left( \begin{array}{cccc}
            a_{11} & \ldots & a_{1n} & 0\\
            \vdots &        & \vdots & \vdots\\
            a_{n1} & \ldots & a_{nn} & 0\\
            0      & \ldots & 0      & 0
  \end{array} \right).
\end{eqnarray}
This is a one-to-one $\ast$-homomorphism, and it is norm preserving.
Define $M_\infty(A)$ as the limit of $M_n(A)$ with respect to these
inclusions. That is, $M_\infty(A)$ is the set of infinite matrices
where almost all $a_{ij}$ are zero. Finally, define the
\emph{stabilization} of $A$ (cf. \cite[6.4]{RLL00} or \cite[1.10]{Wegge-Olsen93}) as the closure
$$
    \dot{A} = \overline{M_\infty(A)}\,.
$$
Here we mean the completion with respect to the norm on $M_\infty(A)$ and the main point is that the inclusions above are all norm-preserving. The result is a \csa{} without unit.

\begin{remark}
There is an equivalent definition of $\dot{A}$ as the tensor product $A \otimes \calk$, where $\calk$ is the \csa{} of all compact operators on a separable infinite-dimensional Hilbert space, and the tensor product is in the sense of \csa{}s.
\end{remark}

\begin{example}
Let $H$ be a separable, infinite-dimensional, Hilbert space.
That is, $H$ has a countable, but not finite, orthonormal basis. It can be shown that
$$
    \dot{\C} = \calk \subset \call(H),
$$
where $\calk$ is the subset of compact operators on $H$. We
have then
$$
    K_j(\C) = K_j(\dot{\C})\,,
$$
where $K_j(-)$ is \csa{} $K$-theory. This is true in general for any \csa{} (Proposition \ref{prop:KtheoryStab} below).\\
On the other hand, the algebraic $K$-theory of $\dot{\C}$ is
$$
    K_j^{\textup{alg}}(\dot{\C}) = \left\{ \begin{array}{cl}
            \Z & \ j \textup{ even,}\\
            0  & \ j \textup{ odd,} \end{array} \right.
$$
which therefore coincides with the \csa{} $K$-theory of \C. This is also true in general (Theorem \ref{thm:KaroubiConj} below). This answer is simple compared with the algebraic $K$-theory of \C{}, where only some partial results are known.
\end{example}

The stabilization of a \csa{} does not change its (\csa{}) $K$-theory.
\begin{proposition}\label{prop:KtheoryStab}
Let $A$ be a \csa{} and write $K_j(-)$ for $K$-theory of \csa{}s.
Then
$$
    K_j(A) = K_j(\dot{A}) \quad j \ge 0\,.
$$
\end{proposition}
The proof is a consequence of the definition of \csa{} $K$-theory: the inclusions (\ref{eqn:Inclusion}) induce isomorphisms in $K$-theory, and the direct limit (in the sense of \csa{}s) commute with the $K$-theory functor (cf. \cite[6.2.11 and 7.1.9]{Wegge-Olsen93}).

\begin{remark}
In the terminology of Corti\~nas' notes\cite{Willie}, Proposition \ref{prop:KtheoryStab} says that the functors $K_0$ and $K_1$ are \emph{$\calk$-stable}.
\end{remark}

M.~Karoubi conjectured that the algebraic $K$-theory of $\dot{A}$ is isomorphic to its \csa{} $K$-theory.
The conjecture was proved by A.~Suslin and M.~Wodzicki.

\begin{theorem}[A.~Suslin and M.~Wodzicki \cite{SW92}]\label{thm:KaroubiConj}
Let $A$ be a \csa. Then
$$
    K_j(\dot{A}) = K^{\textup{alg}}_j(\dot{A}) \quad j \ge 0\,,
$$
where the left-hand side is \csa{} $K$-theory and the right-hand
side is (Quillen's) algebraic $K$-theory of rings.
\end{theorem}
A proof can be found in Corti\~nas' notes \cite[Thm.~7.1.3]{Willie}. In these notes Corti\~nas elaborates the isomorphism above into a long exact sequence which involves cyclic homology.

Theorem \ref{thm:KaroubiConj} is the unity of $K$-theory: It says that \csa{}
$K$-theory is a pleasant subdiscipline of algebraic $K$-theory in
which Bott periodicity is valid and certain basic examples are easy
to calculate.


 \bibliographystyle{amsplain}
 \bibliography{bibliosedano}

\providecommand{\bysame}{\leavevmode\hbox to3em{\hrulefill}\thinspace}
\providecommand{\MR}{\relax\ifhmode\unskip\space\fi MR }
\providecommand{\MRhref}[2]{%
  \href{http://www.ams.org/mathscinet-getitem?mr=#1}{#2}
}
\providecommand{\href}[2]{#2}
\begin{thebibliography}{10}

\bibitem{Arveson76}
William Arveson, \emph{An invitation to {$C\sp*$}-algebras}, Springer-Verlag,
  New York, 1976, Graduate Texts in Mathematics, No. 39. \MR{MR0512360 (58
  \#23621)}

\bibitem{Atiyah70}
M.~F. Atiyah, \emph{Global theory of elliptic operators}, Proc. {I}nternat.
  {C}onf. on {F}unctional {A}nalysis and {R}elated {T}opics ({T}okyo, 1969),
  Univ. of Tokyo Press, Tokyo, 1970, pp.~21--30. \MR{MR0266247 (42 \#1154)}

\bibitem{AH59}
M.~F. Atiyah and F.~Hirzebruch, \emph{Riemann-{R}och theorems for
  differentiable manifolds}, Bull. Amer. Math. Soc. \textbf{65} (1959),
  276--281. \MR{MR0110106 (22 \#989)}

\bibitem{sedanoproceedings}
P.~Baum, G.~Corti\~nas, C.~Mazza, R.~Meyer, M.~Schlichting, and B.~Toen,
  \emph{Lecture {N}otes of the {S}edano {W}inter {S}chool on {$K$}-theory}, to
  appear in Springer Lecture Notes.

\bibitem{BC88}
Paul Baum and Alain Connes, \emph{Chern character for discrete groups}, A
  f\^ete of topology, Academic Press, Boston, MA, 1988, pp.~163--232.
  \MR{MR928402 (90e:58149)}

\bibitem{BCH94}
Paul Baum, Alain Connes, and Nigel Higson, \emph{Classifying space for proper
  actions and {$K$}-theory of group {$C\sp \ast$}-algebras}, $C\sp
  \ast$-algebras: 1943--1993 (San Antonio, TX, 1993), Contemp. Math., vol. 167,
  Amer. Math. Soc., Providence, RI, 1994, pp.~240--291. \MR{MR1292018
  (96c:46070)}

\bibitem{BHS}
Paul Baum, Nigel Higson, and Thomas Schick, \emph{A geometric description of
  equivariant {$K$}-homology}, to appear in the Conference Proceedings of the
  Non Commutative Geometry Conference in honor of Alain Connes, Paris, 2007.

\bibitem{Biller04}
Harald Biller, \emph{Characterizations of proper actions}, Math. Proc.
  Cambridge Philos. Soc. \textbf{136} (2004), no.~2, 429--439. \MR{MR2040583
  (2004k:57043)}

\bibitem{Blackadar98}
Bruce Blackadar, \emph{{$K$}-theory for operator algebras}, second ed.,
  Mathematical Sciences Research Institute Publications, vol.~5, Cambridge
  University Press, Cambridge, 1998. \MR{MR1656031 (99g:46104)}

\bibitem{Bott59}
Raoul Bott, \emph{The stable homotopy of the classical groups}, Ann. of Math.
  (2) \textbf{70} (1959), 313--337. \MR{MR0110104 (22 \#987)}

\bibitem{Brown82}
Kenneth~S. Brown, \emph{Cohomology of groups}, Graduate Texts in Mathematics,
  vol.~87, Springer-Verlag, New York, 1982. \MR{MR672956 (83k:20002)}

\bibitem{CEM01}
J{\'e}r{\^o}me Chabert, Siegfried Echterhoff, and Ralf Meyer, \emph{Deux
  remarques sur l'application de {B}aum-{C}onnes}, C. R. Acad. Sci. Paris
  S\'er. I Math. \textbf{332} (2001), no.~7, 607--610. \MR{MR1841893
  (2002k:19004)}

\bibitem{Willie}
Guillermo Corti\~nas, \emph{Algebraic v.~topological {$K$}-theory: a friendly
  match}, 2009, arXiv:0903.3983v1.

\bibitem{DSV03}
Giuliana Davidoff, Peter Sarnak, and Alain Valette, \emph{Elementary number
  theory, group theory, and {R}amanujan graphs}, London Mathematical Society
  Student Texts, vol.~55, Cambridge University Press, Cambridge, 2003.
  \MR{MR1989434 (2004f:11001)}

\bibitem{DL98}
James~F. Davis and Wolfgang L{\"u}ck, \emph{Spaces over a category and assembly
  maps in isomorphism conjectures in {$K$}- and {$L$}-theory}, $K$-Theory
  \textbf{15} (1998), no.~3, 201--252. \MR{MR1659969 (99m:55004)}

\bibitem{Dixmier77}
Jacques Dixmier, \emph{{$C\sp*$}-algebras}, North-Holland Publishing Co.,
  Amsterdam, 1977, Translated from the French by Francis Jellett, North-Holland
  Mathematical Library, Vol. 15. \MR{MR0458185 (56 \#16388)}

\bibitem{EM09}
Heath Emerson and Ralf Meyer, \emph{Equivariant representable {K}-theory}, J.
  Topol. \textbf{2} (2009), no.~2, 123--156. \MR{MR2499440}

\bibitem{FG05}
Eric~M. Friedlander and Daniel~R. Grayson (eds.), \emph{Handbook of
  {$K$}-theory. {V}ol. 1, 2}, Springer-Verlag, Berlin, 2005. \MR{MR2182598
  (2006e:19001)}

\bibitem{GHW05}
Erik Guentner, Nigel Higson, and Shmuel Weinberger, \emph{The {N}ovikov
  conjecture for linear groups}, Publ. Math. Inst. Hautes \'Etudes Sci. (2005),
  no.~101, 243--268. \MR{MR2217050 (2007c:19007)}

\bibitem{HP04}
Ian Hambleton and Erik~K. Pedersen, \emph{Identifying assembly maps in {$K$}-
  and {$L$}-theory}, Math. Ann. \textbf{328} (2004), no.~1-2, 27--57.
  \MR{MR2030369 (2004j:19001)}

\bibitem{HLS02}
N.~Higson, V.~Lafforgue, and G.~Skandalis, \emph{Counterexamples to the
  {B}aum-{C}onnes conjecture}, Geom. Funct. Anal. \textbf{12} (2002), no.~2,
  330--354. \MR{MR1911663 (2003g:19007)}

\bibitem{Higson90}
Nigel Higson, \emph{A primer on {$KK$}-theory}, Operator theory: operator
  algebras and applications, Part 1 (Durham, NH, 1988), Proc. Sympos. Pure
  Math., vol.~51, Amer. Math. Soc., Providence, RI, 1990, pp.~239--283.
  \MR{MR1077390 (92g:19005)}

\bibitem{HK97}
Nigel Higson and Gennadi Kasparov, \emph{Operator {$K$}-theory for groups which
  act properly and isometrically on {H}ilbert space}, Electron. Res. Announc.
  Amer. Math. Soc. \textbf{3} (1997), 131--142 (electronic). \MR{MR1487204
  (99e:46090)}

\bibitem{HK01}
\bysame, \emph{{$E$}-theory and {$KK$}-theory for groups which act properly and
  isometrically on {H}ilbert space}, Invent. Math. \textbf{144} (2001), no.~1,
  23--74. \MR{MR1821144 (2002k:19005)}

\bibitem{Julg02}
Pierre Julg, \emph{La conjecture de {B}aum-{C}onnes \`a coefficients pour le
  groupe {${\rm Sp}(n,1)$}}, C. R. Math. Acad. Sci. Paris \textbf{334} (2002),
  no.~7, 533--538. \MR{MR1903759 (2003d:19007)}

\bibitem{Kasparov80}
G.~G. Kasparov, \emph{The operator {$K$}-functor and extensions of {$C\sp{\ast}
  $}-algebras}, Izv. Akad. Nauk SSSR Ser. Mat. \textbf{44} (1980), no.~3,
  571--636, 719. \MR{MR582160 (81m:58075)}

\bibitem{Kasparov88}
\bysame, \emph{Equivariant {$KK$}-theory and the {N}ovikov conjecture}, Invent.
  Math. \textbf{91} (1988), no.~1, 147--201. \MR{MR918241 (88j:58123)}

\bibitem{Kasparov95}
\bysame, \emph{{$K$}-theory, group {$C\sp *$}-algebras, and higher signatures
  (conspectus)}, Novikov conjectures, index theorems and rigidity, {V}ol.\ 1
  ({O}berwolfach, 1993), London Math. Soc. Lecture Note Ser., vol. 226,
  Cambridge Univ. Press, Cambridge, 1995, pp.~101--146. \MR{MR1388299
  (97j:58153)}

\bibitem{KS91}
G.~G. Kasparov and G.~Skandalis, \emph{Groups acting on buildings, operator
  {$K$}-theory, and {N}ovikov's conjecture}, $K$-Theory \textbf{4} (1991),
  no.~4, 303--337. \MR{MR1115824 (92h:19009)}

\bibitem{Lafforgue02}
Vincent Lafforgue, \emph{{$K$}-th\'eorie bivariante pour les alg\`ebres de
  {B}anach et conjecture de {B}aum-{C}onnes}, Invent. Math. \textbf{149}
  (2002), no.~1, 1--95. \MR{MR1914617 (2003d:19008)}

\bibitem{Lance95}
E.~C. Lance, \emph{Hilbert {$C\sp *$}-modules}, London Mathematical Society
  Lecture Note Series, vol. 210, Cambridge University Press, Cambridge, 1995, A
  toolkit for operator algebraists. \MR{MR1325694 (96k:46100)}

\bibitem{Luck89}
Wolfgang L{\"u}ck, \emph{Transformation groups and algebraic {$K$}-theory},
  Lecture Notes in Mathematics, vol. 1408, Springer-Verlag, Berlin, 1989,
  Mathematica Gottingensis. \MR{MR1027600 (91g:57036)}

\bibitem{Luck02}
\bysame, \emph{Chern characters for proper equivariant homology theories and
  applications to {$K$}- and {$L$}-theory}, J. Reine Angew. Math. \textbf{543}
  (2002), 193--234.

\bibitem{LuckSurvey}
\bysame, \emph{Survey on classifying spaces for families of subgroups},
  Infinite groups: geometric, combinatorial and dynamical aspects, Progr.
  Math., vol. 248, Birkh\"auser, Basel, 2005, pp.~269--322. \MR{MR2195456
  (2006m:55036)}

\bibitem{LO01}
Wolfgang L{\"u}ck and Bob Oliver, \emph{The completion theorem in {$K$}-theory
  for proper actions of a discrete group}, Topology \textbf{40} (2001), no.~3,
  585--616. \MR{MR1838997 (2002f:19010)}

\bibitem{LR05}
Wolfgang L{\"u}ck and Holger Reich, \emph{The {B}aum-{C}onnes and the
  {F}arrell-{J}ones conjectures in {$K$}- and {$L$}-theory}, Handbook of
  {$K$}-theory. {V}ol. 1, 2, Springer, Berlin, 2005, pp.~703--842.
  \MR{MR2181833 (2006k:19012)}

\bibitem{MOP08}
Michel Matthey, Herv{\'e} Oyono-Oyono, and Wolfgang Pitsch, \emph{Homotopy
  invariance of higher signatures and 3-manifold groups}, Bull. Soc. Math.
  France \textbf{136} (2008), no.~1, 1--25. \MR{MR2415334}

\bibitem{Meyer}
Ralf Meyer, \emph{Universal {C}oefficient {T}heorems and assembly maps in
  {$KK$}-theory}, to appear in \cite{sedanoproceedings}.

\bibitem{MY02}
Igor Mineyev and Guoliang Yu, \emph{The {B}aum-{C}onnes conjecture for
  hyperbolic groups}, Invent. Math. \textbf{149} (2002), no.~1, 97--122.
  \MR{MR1914618 (2003f:20072)}

\bibitem{Mislin03}
Guido Mislin, \emph{Equivariant {$K$}-homology of the classifying space for
  proper actions}, Proper group actions and the Baum-Connes conjecture, Adv.
  Courses Math. CRM Barcelona, Birkh\"auser, Basel, 2003, pp.~1--78.
  \MR{MR2027169 (2005e:19008)}

\bibitem{Murphy90}
Gerard~J. Murphy, \emph{{$C\sp *$}-algebras and operator theory}, Academic
  Press Inc., Boston, MA, 1990. \MR{MR1074574 (91m:46084)}

\bibitem{Oyono-Oyono01}
Herv{\'e} Oyono-Oyono, \emph{Baum-{C}onnes conjecture and group actions on
  trees}, $K$-Theory \textbf{24} (2001), no.~2, 115--134. \MR{MR1869625
  (2002m:19004)}

\bibitem{Pedersen79}
Gert~K. Pedersen, \emph{{$C\sp{\ast} $}-algebras and their automorphism
  groups}, London Mathematical Society Monographs, vol.~14, Academic Press Inc.
  [Harcourt Brace Jovanovich Publishers], London, 1979. \MR{MR548006
  (81e:46037)}

\bibitem{RLL00}
M.~R{\o}rdam, F.~Larsen, and N.~Laustsen, \emph{An introduction to {$K$}-theory
  for {$C\sp *$}-algebras}, London Mathematical Society Student Texts, vol.~49,
  Cambridge University Press, Cambridge, 2000. \MR{MR1783408 (2001g:46001)}

\bibitem{Rosenberg94}
Jonathan Rosenberg, \emph{Algebraic {$K$}-theory and its applications},
  Graduate Texts in Mathematics, vol. 147, Springer-Verlag, New York, 1994.
  \MR{MR1282290 (95e:19001)}

\bibitem{Sanchez-Garcia08}
Rub{\'e}n S{\'a}nchez-Garc{\'{\i}}a, \emph{Bredon homology and equivariant
  {$K$}-homology of {${\rm SL}(3,{\Bbb Z})$}}, J. Pure Appl. Algebra
  \textbf{212} (2008), no.~5, 1046--1059. \MR{MR2387584}

\bibitem{Schick07}
Thomas Schick, \emph{Finite group extensions and the {B}aum-{C}onnes
  conjecture}, Geom. Topol. \textbf{11} (2007), 1767--1775. \MR{MR2350467}

\bibitem{SerreTrees03}
Jean-Pierre Serre, \emph{Trees}, Springer Monographs in Mathematics,
  Springer-Verlag, Berlin, 2003, Translated from the French original by John
  Stillwell, Corrected 2nd printing of the 1980 English translation.
  \MR{MR1954121 (2003m:20032)}

\bibitem{Srinivas96}
V.~Srinivas, \emph{Algebraic {$K$}-theory}, second ed., Progress in
  Mathematics, vol.~90, Birkh\"auser Boston Inc., Boston, MA, 1996.
  \MR{MR1382659 (97c:19001)}

\bibitem{SW92}
Andrei~A. Suslin and Mariusz Wodzicki, \emph{Excision in algebraic
  {$K$}-theory}, Ann. of Math. (2) \textbf{136} (1992), no.~1, 51--122.
  \MR{MR1173926 (93i:19006)}

\bibitem{Wegge-Olsen93}
N.~E. Wegge-Olsen, \emph{{$K$}-theory and {$C\sp *$}-algebras}, Oxford Science
  Publications, The Clarendon Press Oxford University Press, New York, 1993, A
  friendly approach. \MR{MR1222415 (95c:46116)}

\bibitem{WeibelK-Book}
Charles~A. Weibel, \emph{The {$K$}-book: {A}n introduction to algebraic
  {$K$}-theory (book in progress)}, Available at
  \verb+http://www.math.rutgers.edu/~weibel/Kbook.html+.

\bibitem{Weibel99}
\bysame, \emph{The development of algebraic {$K$}-theory before 1980}, Algebra,
  {$K$}-theory, groups, and education ({N}ew {Y}ork, 1997), Contemp. Math.,
  vol. 243, Amer. Math. Soc., Providence, RI, 1999, pp.~211--238. \MR{MR1732049
  (2000m:19001)}

\end{thebibliography}

%
%
\end{document}